\documentclass[12pt]{amsart}
\usepackage[colorlinks=true]{hyperref}
\usepackage[utf8]{inputenc}
\usepackage{amsmath,amsthm,amssymb,url,mathdots,verbatim}
\usepackage[margin = 2.5cm]{geometry}
\usepackage{tikz}
\usepackage{ytableau}
\usepackage{color}
\usepackage{MnSymbol}
\usepackage{graphicx}
\usepackage{scalerel}
\usepackage{enumitem}
\usepackage[normalem]{ulem}
\newtheorem{theorem}{Theorem}

\newtheorem{definition}[theorem]{Definition}

\theoremstyle{definition}

\newtheorem{problem}[theorem]{Problem}

\renewcommand{\u}[1]{\underline{#1}}
\newcommand{\p}[1]{\mathcal #1}
\newcommand{\q}[1]{\mathbf #1}

\DeclareMathOperator{\GT}{\underline{GT}}
\DeclareMathOperator{\MT}{\underline{MT}}
\DeclareMathOperator{\SGT}{\underline{SGT}}
\DeclareMathOperator{\AP}{\underline{AP}}
\DeclareMathOperator{\AR}{\underline{AR}}
\DeclareMathOperator{\DPP}{DPP}
\DeclareMathOperator{\ASM}{ASM}
\DeclareMathOperator{\B}{B}
\DeclareMathOperator{\sign}{sign}

\DeclareMathOperator{\rot}{rot}
\newcommand{\Z}{{\mathbb Z}}
\newcommand{\N}{{\mathbb N}}
\newcommand{\seswarrow}{\swarrow \!\!\!\!\!\;\!\! \searrow}
\newcommand{\nenwarrow}{\nwarrow \!\!\!\!\!\;\!\! \nearrow}
\usepackage{enumitem}
\newcommand{\si}[2]{\u{[#1,#2]}}
\newcommand{\wt}{\widetilde}

\title[A bijective proof of the ASM theorem, Part II]{A bijective proof of the ASM theorem\\ Part II: ASM enumeration and ASM-DPP relation}

\author{Ilse Fischer}
\address{Faculty of Mathematics, University of Vienna, Austria}
\urladdr{https://www.mat.univie.ac.at/~ifischer/}
\author{Matja\v z Konvalinka}%
\address{Faculty of Mathematics and Physics, University of Ljubljana, and Institute of Mathematics, Physics and Mechanics, Slovenia}
\urladdr{http://www.fmf.uni-lj.si/~konvalinka/}
\thanks{The first author acknowledges the financial support from the Austrian Science Foundation FWF, SFB grant F50, and the second author acknowledges the financial support from the Slovenian Research Agency (research core funding No. P1-0294).}

\date{\today}
\begin{document}

\begin{abstract}
This paper is the second in a series of planned papers which provide first bijective proofs of alternating sign matrix results. Based on the main result from the first paper, we construct a bijective proof of the enumeration formula for alternating sign matrices and of the fact that alternating sign matrices are equinumerous with descending plane partitions. We are also able to refine these bijections by including the position of the unique $1$ in the top row of the matrix. Our constructions rely on signed sets and related notions. The starting point for these constructions were known ``computational'' proofs, but the combinatorial point of view led to several drastic modifications.  We also provide computer code where all of our constructions have been implemented.
\end{abstract}

\maketitle

\section{Introduction}

{Bijective combinatorics} is an area of mathematics merely concerned with proving that two finite sets have the same cardinality by constructing explicit bijections between the two sets. Such proofs can be particularly satisfying, especially when the construction is elegant, but often they also reveal many more details about the relation between the two sets than just equinumerosity. A famous set of problems that has resisted numerous attempts to find bijective proofs for more than 35 years now is concerned with \emph{alternating sign matrices} (ASMs) and their relations to certain classes of plane partitions. In his collection of bijective proof problems (which is available from his webpage) Stanley says the following about the problem of finding these bijections: ``\emph{This is one of the most intriguing open problems in the area of bijective proofs.}''
 %In Krattenthaler's survey on plane partitions \cite{krattsurvey} he expresses his opinion by saying: ``\emph{The greatest, still unsolved, mystery concerns the question of what plane partitions have to do with alternating sign matrices.}''
 The current paper is now the second in a planned series of papers that seek to give first bijective proofs of ASM results.

Alternating sign matrices were introduced by Robbins and Rumsey in \cite{lambda} as square matrices with entries in $\{0,1,-1\}$ such that, in each row and each column, the non-zero entries alternate and sum to $1$. They conjectured that the number of $n \times n$ ASMs is given by
$$
\prod_{j=0}^{n-1} \frac{(3j+1)!}{(n+j)!},
$$
which was then proved by Zeilberger in \cite{Zei96a} and a bit later also by
Kuperberg \cite{Kup96}.  For a more extensive account on the history of ASMs, see the first part \cite{PartI}. In this second part we will use the main result from the first paper \cite{PartI} to give a bijective proof of the enumeration formula for $n \times n$ ASMs as well as of the enumeration formula for $n \times n$ ASMs with fixed position of the unique $1$ in the top row. Secondly, we will provide a bijective proof of the fact that there is the same number of \emph{descending plane partitions} (DPPs) with parts less than or equal to $n$ as there is of $n \times n$ alternating sign matrices, and also here we will be able to refine the bijection to include one parameter.

In order to define DPPs, recall that a \emph{strict partition} is a partition $\lambda=(\lambda_1,\ldots,\lambda_l)$ with distinct parts, i.e., $\lambda_1 > \lambda_2 > \ldots > \lambda_l >0$. The shifted Young diagram of shape $\lambda$ is an array of cells with $\lambda_i$ cells in row $i$ and where each row is indented by one cell to the right with respect to the previous row. The shifted Young diagram of shape $(5,3,2)$ is displayed next.
$$
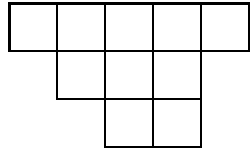
\ydiagram{0+5,1+3,2+2}
$$
A \emph{column strict shifted plane partition} (CSSPP) is a filling of a shifted Young diagram with positive integers such that rows decrease weakly and columns decrease strictly. A descending plane partition (DPP) is a CSSPP such that the first part in each row is greater than the length of its row and less than or equal to the length of the previous row, see \cite{And79}.  An example is given next.
$$
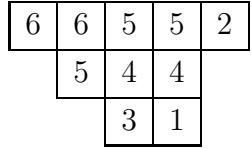
\begin{ytableau}
6 & 6 & 5 & 5 & 2 \\
\none & 5 & 4 & 4  \\
\none & \none & 3 & 1
\end{ytableau}
$$
%It is useful to note that descending plane partitions (DPP) are in easy bijective correspondence with CSSPPs of class $2$, where a CSSPP is of class $k$ if the first part of each row exceeds the length of the row by precisely $k$, see \cite{MRR87}. The CSSPP of class $2$ displayed next corresponds to the DPP above.
%$$
%\begin{ytableau}
%7 & 7 & 6 & 6 & 3 \\
%\none & 6 & 5 & 5 & 1 \\
%\none & \none & 4 & 2
%\end{ytableau}
%$$
%Indeed, to obtain the DPP from the CSSPP of class $2$ one has to subtract $1$ from each part and delete all the zeros then.

\medskip

\subsection*{The ASM enumeration}
The bijection that underlies the bijective proof of the enumeration formula of ASMs as well as the one of the refined enumeration formula involves the following sets:

\begin{itemize}
\item Let $\ASM_{n}$ denote the set of ASMs of size $n \times n$, and, for $1 \le i \le n$, let $\ASM_{n,i}$ denote the subset of $\ASM_n$ of matrices that have the unique $1$ in the first row in column $i$. There is an obvious bijection $\ASM_{n,1} \to \ASM_{n-1}$ which consists of deleting the first row and first column.
\item Let $\B_n$ denote the set of $(2n-1)$-subsets of $[3n-2]=\{1,2,\ldots,3n-2\}$ and, for $1 \le i \le n$, let $\B_{n,i}$ denote the subset of $\B_n$ of those subsets whose median is $n+i-1$. Clearly, $|\B_n| = \binom{3n-2}{2n-1}$ and $|\B_{n,i}| = \binom{n+i-2}{n-1} \binom{2n-i-1}{n-1}$.
\item Let $\DPP_n$ denote the set of descending plane partitions with parts no greater than $n$; let $\DPP_{n,i}$ denote the subset of descending plane partitions with $i-1$ occurrences of $n$. We clearly have $\DPP_{n,1} = \DPP_{n-1}$.
\end{itemize}

\medskip

One main achievement of this paper is the construction of the following bijections.
To emphasize that we are not merely interested in the fact that two signed sets have the same size, but want to use the constructed signed bijection later on, we will be using a convention that is slightly unorthodox in our field. Instead of listing our results as lemmas and theorems with their corresponding proofs, we will be using the Problem--Construction terminology as in \cite{PartI}. See for instance \cite{Voevodsky} and \cite{Bauer}.

\begin{problem} \label{prob:main}
 Given $n \in \N$, $1 \leq i \leq n$, construct a bijection
 $$\DPP_{n-1} \times \B_{n,1} \times \ASM_{n,i} \longrightarrow \DPP_{n-1} \times \ASM_{n,1} \times \B_{n,i}.$$
\end{problem}

Assume that we have constructed such bijections. Then we also have a bijection
\begin{multline*}
\DPP_{n-1} \times \B_{n,1} \times \ASM_n = \bigcup_i \left( \DPP_{n-1} \times  \B_{n,1} \times \ASM_{n,i} \right) \\
\longrightarrow \bigcup_i \left(\DPP_{n-1} \times \ASM_{n,1} \times \B_{n,i} \right) = \DPP_{n-1} \times \ASM_{n,1} \times \B_n \longrightarrow \DPP_{n-1} \times \ASM_{n-1} \times \B_n
\end{multline*}
for every $n$. But by induction, that gives a bijection
$$\DPP_{0} \times \dots \times \DPP_{n-1} \times \B_{1,1} \times \cdots \times \B_{n,1} \times \ASM_n \longrightarrow \DPP_{0} \times \dots \times \DPP_{n-1} \times \B_1 \times \cdots \times \B_n,$$
which proves the ASM theorem
$$|\ASM_n| = \frac{\prod_{j=1}^n |\B_j|}{\prod_{j=1}^n |\B_{j,1}|} = \frac{\prod_{j=1}^n \binom{3j-2}{2j-1}}{\prod_{j=1}^n \binom{2j-2}{j-1}} = \prod_{j=0}^{n-1} \frac{(3j+1)!}{(n+j)!}$$
and also the refined ASM theorem
$$|\ASM_{n,i}| = \frac{|\ASM_{n-1}| \cdot |\B_{n,i}|}{|\B_{n,1}|} = \frac{\binom{n+i-2}{n-1} \binom{2n-i-1}{n-1}}{\binom{3n-2}{2n-1}} \prod_{j=0}^{n-1} \frac{(3j+1)!}{(n+j)!},$$
using the fact that $\DPP_i$ is non-empty (as it contains the empty DPP).

The bijections from Problem~\ref{prob:main} have been implemented in Python and the code is available at \url{https://www.fmf.uni-lj.si/~konvalinka/asmcode.html}. In fact, the bijection depends on an integer parameter $x$ and we provide the case $n=3$, $i=2$, $x=0$ next.

\smallskip

\begin{center}
\scalebox{0.7}{$\begin{array}{llllllll}\left( \emptyset, 12345,
\begin{smallmatrix} 0 & 1 & 0 \\
1 & 0 & 0 \\
0 & 0 & 1 \\ \end{smallmatrix} \right) & \leftrightarrow & \left( \emptyset,
 \begin{smallmatrix} 1 & 0 & 0 \\
0 & 1 & 0 \\
0 & 0 & 1 \\ \end{smallmatrix},
23457\right) \qquad \left( \emptyset, 12345,
\begin{smallmatrix} 0 & 1 & 0 \\
1 & -1 & 1 \\
0 & 1 & 0 \\ \end{smallmatrix} \right) & \leftrightarrow & \left( \emptyset,
 \begin{smallmatrix} 1 & 0 & 0 \\
0 & 0 & 1 \\
0 & 1 & 0 \\ \end{smallmatrix},
23456\right)  \qquad \left( \emptyset, 12345,
\begin{smallmatrix} 0 & 1 & 0 \\
0 & 0 & 1 \\
1 & 0 & 0 \\ \end{smallmatrix} \right) & \leftrightarrow & \left( \emptyset,
 \begin{smallmatrix} 1 & 0 & 0 \\
0 & 1 & 0 \\
0 & 0 & 1 \\ \end{smallmatrix},
23456\right) \vspace{1mm} \\
\left( \emptyset, 12346,
\begin{smallmatrix} 0 & 1 & 0 \\
1 & 0 & 0 \\
0 & 0 & 1 \\ \end{smallmatrix} \right) & \leftrightarrow & \left( \emptyset,
 \begin{smallmatrix} 1 & 0 & 0 \\
0 & 1 & 0 \\
0 & 0 & 1 \\ \end{smallmatrix},
13457\right) \qquad \left( \emptyset, 12346,
\begin{smallmatrix} 0 & 1 & 0 \\
1 & -1 & 1 \\
0 & 1 & 0 \\ \end{smallmatrix} \right) & \leftrightarrow & \left( \emptyset,
 \begin{smallmatrix} 1 & 0 & 0 \\
0 & 0 & 1 \\
0 & 1 & 0 \\ \end{smallmatrix},
13456\right) \qquad \left( \emptyset, 12346,
\begin{smallmatrix} 0 & 1 & 0 \\
0 & 0 & 1 \\
1 & 0 & 0 \\ \end{smallmatrix} \right) & \leftrightarrow & \left( \emptyset,
 \begin{smallmatrix} 1 & 0 & 0 \\
0 & 1 & 0 \\
0 & 0 & 1 \\ \end{smallmatrix},
13456\right) \vspace{1mm} \\
\left( \emptyset, 12347,
\begin{smallmatrix} 0 & 1 & 0 \\
1 & 0 & 0 \\
0 & 0 & 1 \\ \end{smallmatrix} \right) & \leftrightarrow & \left( \emptyset,
 \begin{smallmatrix} 1 & 0 & 0 \\
0 & 1 & 0 \\
0 & 0 & 1 \\ \end{smallmatrix},
12457\right) \qquad \left( \emptyset, 12347,
\begin{smallmatrix} 0 & 1 & 0 \\
1 & -1 & 1 \\
0 & 1 & 0 \\ \end{smallmatrix} \right) & \leftrightarrow & \left( \emptyset,
 \begin{smallmatrix} 1 & 0 & 0 \\
0 & 0 & 1 \\
0 & 1 & 0 \\ \end{smallmatrix},
12456\right) \qquad \left( \emptyset, 12347,
\begin{smallmatrix} 0 & 1 & 0 \\
0 & 0 & 1 \\
1 & 0 & 0 \\ \end{smallmatrix} \right) & \leftrightarrow & \left( \emptyset,
 \begin{smallmatrix} 1 & 0 & 0 \\
0 & 1 & 0 \\
0 & 0 & 1 \\ \end{smallmatrix},
12456\right) \vspace{1mm} \\
 \left( \emptyset, 12356,
\begin{smallmatrix} 0 & 1 & 0 \\
1 & 0 & 0 \\
0 & 0 & 1 \\ \end{smallmatrix} \right) & \leftrightarrow & \left( 2,
 \begin{smallmatrix} 1 & 0 & 0 \\
0 & 1 & 0 \\
0 & 0 & 1 \\ \end{smallmatrix},
13456\right) \qquad \left( \emptyset, 12356,
\begin{smallmatrix} 0 & 1 & 0 \\
1 & -1 & 1 \\
0 & 1 & 0 \\ \end{smallmatrix} \right) & \leftrightarrow & \left( 2,
 \begin{smallmatrix} 1 & 0 & 0 \\
0 & 0 & 1 \\
0 & 1 & 0 \\ \end{smallmatrix},
12456\right) \qquad \left( \emptyset, 12356,
\begin{smallmatrix} 0 & 1 & 0 \\
0 & 0 & 1 \\
1 & 0 & 0 \\ \end{smallmatrix} \right) & \leftrightarrow & \left( 2,
 \begin{smallmatrix} 1 & 0 & 0 \\
0 & 1 & 0 \\
0 & 0 & 1 \\ \end{smallmatrix},
12456\right) \vspace{1mm} \\ \left( \emptyset, 12357,
\begin{smallmatrix} 0 & 1 & 0 \\
1 & 0 & 0 \\
0 & 0 & 1 \\ \end{smallmatrix} \right) & \leftrightarrow & \left( 2,
 \begin{smallmatrix} 1 & 0 & 0 \\
0 & 1 & 0 \\
0 & 0 & 1 \\ \end{smallmatrix},
13457\right) \qquad \left( \emptyset, 12357,
\begin{smallmatrix} 0 & 1 & 0 \\
1 & -1 & 1 \\
0 & 1 & 0 \\ \end{smallmatrix} \right) & \leftrightarrow & \left( 2,
 \begin{smallmatrix} 1 & 0 & 0 \\
0 & 0 & 1 \\
0 & 1 & 0 \\ \end{smallmatrix},
12457\right) \qquad \left( \emptyset, 12357,
\begin{smallmatrix} 0 & 1 & 0 \\
0 & 0 & 1 \\
1 & 0 & 0 \\ \end{smallmatrix} \right) & \leftrightarrow & \left( 2,
 \begin{smallmatrix} 1 & 0 & 0 \\
0 & 1 & 0 \\
0 & 0 & 1 \\ \end{smallmatrix},
12457\right) \vspace{1mm} \\
 \left( \emptyset, 12367,
\begin{smallmatrix} 0 & 1 & 0 \\
1 & 0 & 0 \\
0 & 0 & 1 \\ \end{smallmatrix} \right) & \leftrightarrow & \left( 2,
 \begin{smallmatrix} 1 & 0 & 0 \\
0 & 1 & 0 \\
0 & 0 & 1 \\ \end{smallmatrix},
13467\right) \qquad \left( \emptyset, 12367,
\begin{smallmatrix} 0 & 1 & 0 \\
1 & -1 & 1 \\
0 & 1 & 0 \\ \end{smallmatrix} \right) & \leftrightarrow & \left( 2,
 \begin{smallmatrix} 1 & 0 & 0 \\
0 & 0 & 1 \\
0 & 1 & 0 \\ \end{smallmatrix},
12467\right) \qquad \left( \emptyset, 12367,
\begin{smallmatrix} 0 & 1 & 0 \\
0 & 0 & 1 \\
1 & 0 & 0 \\ \end{smallmatrix} \right) & \leftrightarrow & \left( 2,
 \begin{smallmatrix} 1 & 0 & 0 \\
0 & 1 & 0 \\
0 & 0 & 1 \\ \end{smallmatrix},
12467\right) \vspace{1mm} \\
\left( 2,
12345,
\begin{smallmatrix} 0 & 1 & 0 \\
1 & 0 & 0 \\
0 & 0 & 1 \\ \end{smallmatrix} \right) & \leftrightarrow & \left( \emptyset,
 \begin{smallmatrix} 1 & 0 & 0 \\
0 & 1 & 0 \\
0 & 0 & 1 \\ \end{smallmatrix},
23467\right) \qquad \left( 2,
12345,
\begin{smallmatrix} 0 & 1 & 0 \\
1 & -1 & 1 \\
0 & 1 & 0 \\ \end{smallmatrix} \right) & \leftrightarrow & \left( \emptyset,
 \begin{smallmatrix} 1 & 0 & 0 \\
0 & 0 & 1 \\
0 & 1 & 0 \\ \end{smallmatrix},
23467\right) \qquad \left( 2,
12345,
\begin{smallmatrix} 0 & 1 & 0 \\
0 & 0 & 1 \\
1 & 0 & 0 \\ \end{smallmatrix} \right) & \leftrightarrow & \left( \emptyset,
 \begin{smallmatrix} 1 & 0 & 0 \\
0 & 0 & 1 \\
0 & 1 & 0 \\ \end{smallmatrix},
23457\right) \vspace{1mm} \\
 \left( 2,
12346,
\begin{smallmatrix} 0 & 1 & 0 \\
1 & 0 & 0 \\
0 & 0 & 1 \\ \end{smallmatrix} \right) & \leftrightarrow & \left( \emptyset,
 \begin{smallmatrix} 1 & 0 & 0 \\
0 & 1 & 0 \\
0 & 0 & 1 \\ \end{smallmatrix},
13467\right) \qquad \left( 2,
12346,
\begin{smallmatrix} 0 & 1 & 0 \\
1 & -1 & 1 \\
0 & 1 & 0 \\ \end{smallmatrix} \right) & \leftrightarrow & \left( \emptyset,
 \begin{smallmatrix} 1 & 0 & 0 \\
0 & 0 & 1 \\
0 & 1 & 0 \\ \end{smallmatrix},
13467\right) \qquad \left( 2,
12346,
\begin{smallmatrix} 0 & 1 & 0 \\
0 & 0 & 1 \\
1 & 0 & 0 \\ \end{smallmatrix} \right) & \leftrightarrow & \left( \emptyset,
 \begin{smallmatrix} 1 & 0 & 0 \\
0 & 0 & 1 \\
0 & 1 & 0 \\ \end{smallmatrix},
13457\right) \vspace{1mm} \\
\left( 2,
12347,
\begin{smallmatrix} 0 & 1 & 0 \\
1 & 0 & 0 \\
0 & 0 & 1 \\ \end{smallmatrix} \right) & \leftrightarrow & \left( \emptyset,
 \begin{smallmatrix} 1 & 0 & 0 \\
0 & 1 & 0 \\
0 & 0 & 1 \\ \end{smallmatrix},
12467\right) \qquad \left( 2,
12347,
\begin{smallmatrix} 0 & 1 & 0 \\
1 & -1 & 1 \\
0 & 1 & 0 \\ \end{smallmatrix} \right) & \leftrightarrow & \left( \emptyset,
 \begin{smallmatrix} 1 & 0 & 0 \\
0 & 0 & 1 \\
0 & 1 & 0 \\ \end{smallmatrix},
12467\right) \qquad \left( 2,
12347,
\begin{smallmatrix} 0 & 1 & 0 \\
0 & 0 & 1 \\
1 & 0 & 0 \\ \end{smallmatrix} \right) & \leftrightarrow & \left( \emptyset,
 \begin{smallmatrix} 1 & 0 & 0 \\
0 & 0 & 1 \\
0 & 1 & 0 \\ \end{smallmatrix},
12457\right) \vspace{1mm} \\
 \left( 2,
12356,
\begin{smallmatrix} 0 & 1 & 0 \\
1 & 0 & 0 \\
0 & 0 & 1 \\ \end{smallmatrix} \right) & \leftrightarrow & \left( 2,
 \begin{smallmatrix} 1 & 0 & 0 \\
0 & 1 & 0 \\
0 & 0 & 1 \\ \end{smallmatrix},
23456\right) \qquad \left( 2,
12356,
\begin{smallmatrix} 0 & 1 & 0 \\
1 & -1 & 1 \\
0 & 1 & 0 \\ \end{smallmatrix} \right) & \leftrightarrow & \left( 2,
 \begin{smallmatrix} 1 & 0 & 0 \\
0 & 0 & 1 \\
0 & 1 & 0 \\ \end{smallmatrix},
23456\right) \qquad \left( 2,
12356,
\begin{smallmatrix} 0 & 1 & 0 \\
0 & 0 & 1 \\
1 & 0 & 0 \\ \end{smallmatrix} \right) & \leftrightarrow & \left( 2,
 \begin{smallmatrix} 1 & 0 & 0 \\
0 & 0 & 1 \\
0 & 1 & 0 \\ \end{smallmatrix},
13456\right) \vspace{1mm} \\
 \left( 2,
12357,
\begin{smallmatrix} 0 & 1 & 0 \\
1 & 0 & 0 \\
0 & 0 & 1 \\ \end{smallmatrix} \right) & \leftrightarrow & \left( 2,
 \begin{smallmatrix} 1 & 0 & 0 \\
0 & 1 & 0 \\
0 & 0 & 1 \\ \end{smallmatrix},
23457\right) \qquad \left( 2,
12357,
\begin{smallmatrix} 0 & 1 & 0 \\
1 & -1 & 1 \\
0 & 1 & 0 \\ \end{smallmatrix} \right) & \leftrightarrow & \left( 2,
 \begin{smallmatrix} 1 & 0 & 0 \\
0 & 0 & 1 \\
0 & 1 & 0 \\ \end{smallmatrix},
23457\right) \qquad \left( 2,
12357,
\begin{smallmatrix} 0 & 1 & 0 \\
0 & 0 & 1 \\
1 & 0 & 0 \\ \end{smallmatrix} \right) & \leftrightarrow & \left( 2,
 \begin{smallmatrix} 1 & 0 & 0 \\
0 & 0 & 1 \\
0 & 1 & 0 \\ \end{smallmatrix},
13457\right) \vspace{1mm} \\ \left( 2,
12367,
\begin{smallmatrix} 0 & 1 & 0 \\
1 & 0 & 0 \\
0 & 0 & 1 \\ \end{smallmatrix} \right) & \leftrightarrow & \left( 2,
 \begin{smallmatrix} 1 & 0 & 0 \\
0 & 1 & 0 \\
0 & 0 & 1 \\ \end{smallmatrix},
23467\right) \qquad \left( 2,
12367,
\begin{smallmatrix} 0 & 1 & 0 \\
1 & -1 & 1 \\
0 & 1 & 0 \\ \end{smallmatrix} \right) & \leftrightarrow & \left( 2,
 \begin{smallmatrix} 1 & 0 & 0 \\
0 & 0 & 1 \\
0 & 1 & 0 \\ \end{smallmatrix},
23467\right) \qquad \left( 2,
12367,
\begin{smallmatrix} 0 & 1 & 0 \\
0 & 0 & 1 \\
1 & 0 & 0 \\ \end{smallmatrix} \right) & \leftrightarrow & \left( 2,
 \begin{smallmatrix} 1 & 0 & 0 \\
0 & 0 & 1 \\
0 & 1 & 0 \\ \end{smallmatrix},
13467\right)\\
\end{array}$}
\end{center}

\subsection*{The ASM-DPP relation}
We also construct an explicit bijection that proves $|\ASM_n| = |\DPP_n|$ and more generally $|\ASM_{n,j}| = |\DPP_{n,j}|$.

\begin{problem} \label{prob:asmtodpp}
 Given $n \in \N$, $1 \leq i \leq n$, construct a bijection
 $$\DPP_{n-1} \times \ASM_{n,i} \longrightarrow \ASM_{n,1} \times \DPP_{n,i}.$$
\end{problem}

Indeed, once this is proven it follows that
$
|\DPP_{n-1}| \cdot |\ASM_{n,i}| = |\ASM_{n,1}| \cdot |\DPP_{n,i}| =  |\ASM_{n-1}| \cdot |\DPP_{n,i}|.
$
By induction, we can assume $|\DPP_{n-1}|=|\ASM_{n-1}|$ and so $|\ASM_{n,i}| = |\DPP_{n,i}|$. Summing this over all $i$ implies $|\DPP_{n}|=|\ASM_{n}|$.

The bijections from Problem~\ref{prob:asmtodpp} have also been implemented, see the above mentioned webpage. Again the bijection depends on an integer parameter $x$ and the case $n=4$, $i=2$, $x=0$ is as follows.

\smallskip

\begin{center}
\scalebox{0.5}{$\begin{array}{cccccccccccc}
\left(\emptyset, \begin{smallmatrix}
 0 &  1 &  0 &  0 \\
 1 &  0 &  0 &  0 \\
 0 &  0 &  1 &  0 \\
 0 &  0 &  0 &  1
\end{smallmatrix}\right)  & \leftrightarrow &
\left(\begin{smallmatrix}
 1 &  0 &  0 &  0 \\
 0 &  1 &  0 &  0 \\
 0 &  0 &  0 &  1 \\
 0 &  0 &  1 &  0
\end{smallmatrix}, \begin{array}{ccc}
4 & 2 & 1
\end{array}\right) & \left(\emptyset, \begin{smallmatrix}
 0 &  1 &  0 &  0 \\
 1 & -1 &  1 &  0 \\
 0 &  1 &  0 &  0 \\
 0 &  0 &  0 &  1
\end{smallmatrix}\right)  & \leftrightarrow &
\left(\begin{smallmatrix}
 1 &  0 &  0 &  0 \\
 0 &  0 &  1 &  0 \\
 0 &  1 & -1 &  1 \\
 0 &  0 &  1 &  0
\end{smallmatrix}, \begin{array}{ccc}
4 & 1 & 1
\end{array}\right) & \left(\emptyset, \begin{smallmatrix}
 0 &  1 &  0 &  0 \\
 0 &  0 &  1 &  0 \\
 1 &  0 &  0 &  0 \\
 0 &  0 &  0 &  1
\end{smallmatrix}\right)  & \leftrightarrow &
\left(\begin{smallmatrix}
 1 &  0 &  0 &  0 \\
 0 &  0 &  0 &  1 \\
 0 &  1 &  0 &  0 \\
 0 &  0 &  1 &  0
\end{smallmatrix}, \begin{array}{cc}
4 & 1
\end{array}\right) &
\left(\emptyset, \begin{smallmatrix}
 0 &  1 &  0 &  0 \\
 1 &  0 &  0 &  0 \\
 0 &  0 &  0 &  1 \\
 0 &  0 &  1 &  0
\end{smallmatrix}\right)  & \leftrightarrow &
\left(\begin{smallmatrix}
 1 &  0 &  0 &  0 \\
 0 &  1 &  0 &  0 \\
 0 &  0 &  1 &  0 \\
 0 &  0 &  0 &  1
\end{smallmatrix}, \begin{array}{cc}
4 & 2
\end{array}\right) \\[0.4cm]
 \left(\emptyset, \begin{smallmatrix}
 0 &  1 &  0 &  0 \\
 1 & -1 &  1 &  0 \\
 0 &  1 & -1 &  1 \\
 0 &  0 &  1 &  0
\end{smallmatrix}\right)  & \leftrightarrow &
\left(\begin{smallmatrix}
 1 &  0 &  0 &  0 \\
 0 &  1 &  0 &  0 \\
 0 &  0 &  1 &  0 \\
 0 &  0 &  0 &  1
\end{smallmatrix}, \begin{array}{cc}
4 & 1
\end{array}\right) & \left(\emptyset, \begin{smallmatrix}
 0 &  1 &  0 &  0 \\
 1 & -1 &  0 &  1 \\
 0 &  1 &  0 &  0 \\
 0 &  0 &  1 &  0
\end{smallmatrix}\right)  & \leftrightarrow &
\left(\begin{smallmatrix}
 1 &  0 &  0 &  0 \\
 0 &  1 &  0 &  0 \\
 0 &  0 &  1 &  0 \\
 0 &  0 &  0 &  1
\end{smallmatrix}, \begin{array}{c}
4
\end{array}\right) &
\left(\emptyset, \begin{smallmatrix}
 0 &  1 &  0 &  0 \\
 0 &  0 &  1 &  0 \\
 1 &  0 & -1 &  1 \\
 0 &  0 &  1 &  0
\end{smallmatrix}\right)  & \leftrightarrow &
\left(\begin{smallmatrix}
 1 &  0 &  0 &  0 \\
 0 &  0 &  1 &  0 \\
 0 &  0 &  0 &  1 \\
 0 &  1 &  0 &  0
\end{smallmatrix}, \begin{array}{c}
4
\end{array}\right) & \left(\emptyset, \begin{smallmatrix}
 0 &  1 &  0 &  0 \\
 0 &  0 &  0 &  1 \\
 1 &  0 &  0 &  0 \\
 0 &  0 &  1 &  0
\end{smallmatrix}\right)  & \leftrightarrow &
\left(\begin{smallmatrix}
 1 &  0 &  0 &  0 \\
 0 &  1 &  0 &  0 \\
 0 &  0 &  0 &  1 \\
 0 &  0 &  1 &  0
\end{smallmatrix}, \begin{array}{ccc}
4 & 1 & 1
\end{array}\right) \\[0.4cm]
 \left(\emptyset, \begin{smallmatrix}
 0 &  1 &  0 &  0 \\
 1 & -1 &  1 &  0 \\
 0 &  0 &  0 &  1 \\
 0 &  1 &  0 &  0
\end{smallmatrix}\right)  & \leftrightarrow &
\left(\begin{smallmatrix}
 1 &  0 &  0 &  0 \\
 0 &  1 &  0 &  0 \\
 0 &  0 &  1 &  0 \\
 0 &  0 &  0 &  1
\end{smallmatrix}, \begin{array}{ccc}
4 & 1 & 1
\end{array}\right) &
\left(\emptyset, \begin{smallmatrix}
 0 &  1 &  0 &  0 \\
 1 & -1 &  0 &  1 \\
 0 &  0 &  1 &  0 \\
 0 &  1 &  0 &  0
\end{smallmatrix}\right)  & \leftrightarrow &
\left(\begin{smallmatrix}
 1 &  0 &  0 &  0 \\
 0 &  0 &  0 &  1 \\
 0 &  0 &  1 &  0 \\
 0 &  1 &  0 &  0
\end{smallmatrix}, \begin{array}{c}
4
\end{array}\right) & \left(\emptyset, \begin{smallmatrix}
 0 &  1 &  0 &  0 \\
 0 &  0 &  1 &  0 \\
 1 & -1 &  0 &  1 \\
 0 &  1 &  0 &  0
\end{smallmatrix}\right)  & \leftrightarrow &
\left(\begin{smallmatrix}
 1 &  0 &  0 &  0 \\
 0 &  1 &  0 &  0 \\
 0 &  0 &  0 &  1 \\
 0 &  0 &  1 &  0
\end{smallmatrix}, \begin{array}{cc}
4 & 2
\end{array}\right) &
 \left(\emptyset, \begin{smallmatrix}
 0 &  1 &  0 &  0 \\
 0 &  0 &  0 &  1 \\
 1 & -1 &  1 &  0 \\
 0 &  1 &  0 &  0
\end{smallmatrix}\right)  & \leftrightarrow &
\left(\begin{smallmatrix}
 1 &  0 &  0 &  0 \\
 0 &  1 &  0 &  0 \\
 0 &  0 &  0 &  1 \\
 0 &  0 &  1 &  0
\end{smallmatrix}, \begin{array}{c}
4
\end{array}\right)\\[0.4cm]
\left(\emptyset, \begin{smallmatrix}
 0 &  1 &  0 &  0 \\
 0 &  0 &  1 &  0 \\
 0 &  0 &  0 &  1 \\
 1 &  0 &  0 &  0
\end{smallmatrix}\right)  & \leftrightarrow &
\left(\begin{smallmatrix}
 1 &  0 &  0 &  0 \\
 0 &  0 &  1 &  0 \\
 0 &  1 & -1 &  1 \\
 0 &  0 &  1 &  0
\end{smallmatrix}, \begin{array}{cc}
4 & 1
\end{array}\right) & \left(\emptyset, \begin{smallmatrix}
 0 &  1 &  0 &  0 \\
 0 &  0 &  0 &  1 \\
 0 &  0 &  1 &  0 \\
 1 &  0 &  0 &  0
\end{smallmatrix}\right)  & \leftrightarrow &
\left(\begin{smallmatrix}
 1 &  0 &  0 &  0 \\
 0 &  1 &  0 &  0 \\
 0 &  0 &  0 &  1 \\
 0 &  0 &  1 &  0
\end{smallmatrix}, \begin{array}{cc}
4 & 1
\end{array}\right) & \left(\begin{array}{c}
2
\end{array}, \begin{smallmatrix}
 0 &  1 &  0 &  0 \\
 1 &  0 &  0 &  0 \\
 0 &  0 &  1 &  0 \\
 0 &  0 &  0 &  1
\end{smallmatrix}\right)  & \leftrightarrow &
\left(\begin{smallmatrix}
 1 &  0 &  0 &  0 \\
 0 &  1 &  0 &  0 \\
 0 &  0 &  1 &  0 \\
 0 &  0 &  0 &  1
\end{smallmatrix}, \begin{array}{cc}
4 & 3
\end{array}\right) &
\left(\begin{array}{c}
2
\end{array}, \begin{smallmatrix}
 0 &  1 &  0 &  0 \\
 1 & -1 &  1 &  0 \\
 0 &  1 &  0 &  0 \\
 0 &  0 &  0 &  1
\end{smallmatrix}\right)  & \leftrightarrow &
\left(\begin{smallmatrix}
 1 &  0 &  0 &  0 \\
 0 &  0 &  1 &  0 \\
 0 &  1 &  0 &  0 \\
 0 &  0 &  0 &  1
\end{smallmatrix}, \begin{array}{cc}
4 & 3
\end{array}\right) \\[0.4cm]
 \left(\begin{array}{c}
2
\end{array}, \begin{smallmatrix}
 0 &  1 &  0 &  0 \\
 0 &  0 &  1 &  0 \\
 1 &  0 &  0 &  0 \\
 0 &  0 &  0 &  1
\end{smallmatrix}\right)  & \leftrightarrow &
\left(\begin{smallmatrix}
 1 &  0 &  0 &  0 \\
 0 &  0 &  1 &  0 \\
 0 &  1 &  0 &  0 \\
 0 &  0 &  0 &  1
\end{smallmatrix}, \begin{array}{cc}
4 & 2
\end{array}\right) & \left(\begin{array}{c}
2
\end{array}, \begin{smallmatrix}
 0 &  1 &  0 &  0 \\
 1 &  0 &  0 &  0 \\
 0 &  0 &  0 &  1 \\
 0 &  0 &  1 &  0
\end{smallmatrix}\right)  & \leftrightarrow &
\left(\begin{smallmatrix}
 1 &  0 &  0 &  0 \\
 0 &  1 &  0 &  0 \\
 0 &  0 &  0 &  1 \\
 0 &  0 &  1 &  0
\end{smallmatrix}, \begin{array}{cc}
4 & 3
\end{array}\right) &
\left(\begin{array}{c}
2
\end{array}, \begin{smallmatrix}
 0 &  1 &  0 &  0 \\
 1 & -1 &  1 &  0 \\
 0 &  1 & -1 &  1 \\
 0 &  0 &  1 &  0
\end{smallmatrix}\right)  & \leftrightarrow &
\left(\begin{smallmatrix}
 1 &  0 &  0 &  0 \\
 0 &  0 &  1 &  0 \\
 0 &  1 & -1 &  1 \\
 0 &  0 &  1 &  0
\end{smallmatrix}, \begin{array}{cc}
4 & 3
\end{array}\right) & \left(\begin{array}{c}
2
\end{array}, \begin{smallmatrix}
 0 &  1 &  0 &  0 \\
 1 & -1 &  0 &  1 \\
 0 &  1 &  0 &  0 \\
 0 &  0 &  1 &  0
\end{smallmatrix}\right)  & \leftrightarrow &
\left(\begin{smallmatrix}
 1 &  0 &  0 &  0 \\
 0 &  0 &  0 &  1 \\
 0 &  1 &  0 &  0 \\
 0 &  0 &  1 &  0
\end{smallmatrix}, \begin{array}{cc}
4 & 3
\end{array}\right) \\[0.4cm]
 \left(\begin{array}{c}
2
\end{array}, \begin{smallmatrix}
 0 &  1 &  0 &  0 \\
 0 &  0 &  1 &  0 \\
 1 &  0 & -1 &  1 \\
 0 &  0 &  1 &  0
\end{smallmatrix}\right)  & \leftrightarrow &
\left(\begin{smallmatrix}
 1 &  0 &  0 &  0 \\
 0 &  0 &  1 &  0 \\
 0 &  1 & -1 &  1 \\
 0 &  0 &  1 &  0
\end{smallmatrix}, \begin{array}{cc}
4 & 2
\end{array}\right) &
\left(\begin{array}{c}
2
\end{array}, \begin{smallmatrix}
 0 &  1 &  0 &  0 \\
 0 &  0 &  0 &  1 \\
 1 &  0 &  0 &  0 \\
 0 &  0 &  1 &  0
\end{smallmatrix}\right)  & \leftrightarrow &
\left(\begin{smallmatrix}
 1 &  0 &  0 &  0 \\
 0 &  0 &  0 &  1 \\
 0 &  1 &  0 &  0 \\
 0 &  0 &  1 &  0
\end{smallmatrix}, \begin{array}{cc}
4 & 2
\end{array}\right) & \left(\begin{array}{c}
2
\end{array}, \begin{smallmatrix}
 0 &  1 &  0 &  0 \\
 1 & -1 &  1 &  0 \\
 0 &  0 &  0 &  1 \\
 0 &  1 &  0 &  0
\end{smallmatrix}\right)  & \leftrightarrow &
\left(\begin{smallmatrix}
 1 &  0 &  0 &  0 \\
 0 &  0 &  1 &  0 \\
 0 &  0 &  0 &  1 \\
 0 &  1 &  0 &  0
\end{smallmatrix}, \begin{array}{cc}
4 & 3
\end{array}\right) & \left(\begin{array}{c}
2
\end{array}, \begin{smallmatrix}
 0 &  1 &  0 &  0 \\
 1 & -1 &  0 &  1 \\
 0 &  0 &  1 &  0 \\
 0 &  1 &  0 &  0
\end{smallmatrix}\right)  & \leftrightarrow &
\left(\begin{smallmatrix}
 1 &  0 &  0 &  0 \\
 0 &  0 &  0 &  1 \\
 0 &  0 &  1 &  0 \\
 0 &  1 &  0 &  0
\end{smallmatrix}, \begin{array}{cc}
4 & 3
\end{array}\right)\\[0.4cm]
\left(\begin{array}{c}
2
\end{array}, \begin{smallmatrix}
 0 &  1 &  0 &  0 \\
 0 &  0 &  1 &  0 \\
 1 & -1 &  0 &  1 \\
 0 &  1 &  0 &  0
\end{smallmatrix}\right)  & \leftrightarrow &
\left(\begin{smallmatrix}
 1 &  0 &  0 &  0 \\
 0 &  0 &  1 &  0 \\
 0 &  0 &  0 &  1 \\
 0 &  1 &  0 &  0
\end{smallmatrix}, \begin{array}{cc}
4 & 2
\end{array}\right) & \left(\begin{array}{c}
2
\end{array}, \begin{smallmatrix}
 0 &  1 &  0 &  0 \\
 0 &  0 &  0 &  1 \\
 1 & -1 &  1 &  0 \\
 0 &  1 &  0 &  0
\end{smallmatrix}\right)  & \leftrightarrow &
\left(\begin{smallmatrix}
 1 &  0 &  0 &  0 \\
 0 &  0 &  0 &  1 \\
 0 &  0 &  1 &  0 \\
 0 &  1 &  0 &  0
\end{smallmatrix}, \begin{array}{cc}
4 & 2
\end{array}\right) & \left(\begin{array}{c}
2
\end{array}, \begin{smallmatrix}
 0 &  1 &  0 &  0 \\
 0 &  0 &  1 &  0 \\
 0 &  0 &  0 &  1 \\
 1 &  0 &  0 &  0
\end{smallmatrix}\right)  & \leftrightarrow &
\left(\begin{smallmatrix}
 1 &  0 &  0 &  0 \\
 0 &  0 &  1 &  0 \\
 0 &  0 &  0 &  1 \\
 0 &  1 &  0 &  0
\end{smallmatrix}, \begin{array}{cc}
4 & 1
\end{array}\right) &
\left(\begin{array}{c}
2
\end{array}, \begin{smallmatrix}
 0 &  1 &  0 &  0 \\
 0 &  0 &  0 &  1 \\
 0 &  0 &  1 &  0 \\
 1 &  0 &  0 &  0
\end{smallmatrix}\right)  & \leftrightarrow &
\left(\begin{smallmatrix}
 1 &  0 &  0 &  0 \\
 0 &  0 &  0 &  1 \\
 0 &  0 &  1 &  0 \\
 0 &  1 &  0 &  0
\end{smallmatrix}, \begin{array}{cc}
4 & 1
\end{array}\right) \\[0.4cm]
 \left(\begin{array}{cc}
3 & 3
\end{array}, \begin{smallmatrix}
 0 &  1 &  0 &  0 \\
 1 &  0 &  0 &  0 \\
 0 &  0 &  1 &  0 \\
 0 &  0 &  0 &  1
\end{smallmatrix}\right)  & \leftrightarrow &
\left(\begin{smallmatrix}
 1 &  0 &  0 &  0 \\
 0 &  1 &  0 &  0 \\
 0 &  0 &  1 &  0 \\
 0 &  0 &  0 &  1
\end{smallmatrix}, \begin{array}{ccc}
4 & 2 & 1
\end{array}\right) & \left(\begin{array}{cc}
3 & 3
\end{array}, \begin{smallmatrix}
 0 &  1 &  0 &  0 \\
 1 & -1 &  1 &  0 \\
 0 &  1 &  0 &  0 \\
 0 &  0 &  0 &  1
\end{smallmatrix}\right)  & \leftrightarrow &
\left(\begin{smallmatrix}
 1 &  0 &  0 &  0 \\
 0 &  1 &  0 &  0 \\
 0 &  0 &  1 &  0 \\
 0 &  0 &  0 &  1
\end{smallmatrix}, \begin{array}{ccc}
4 & 3 & 1 \\
  & 2
\end{array}\right) &
\left(\begin{array}{cc}
3 & 3
\end{array}, \begin{smallmatrix}
 0 &  1 &  0 &  0 \\
 0 &  0 &  1 &  0 \\
 1 &  0 &  0 &  0 \\
 0 &  0 &  0 &  1
\end{smallmatrix}\right)  & \leftrightarrow &
\left(\begin{smallmatrix}
 1 &  0 &  0 &  0 \\
 0 &  0 &  0 &  1 \\
 0 &  1 &  0 &  0 \\
 0 &  0 &  1 &  0
\end{smallmatrix}, \begin{array}{ccc}
4 & 3 & 2 \\
  & 2
\end{array}\right) & \left(\begin{array}{cc}
3 & 3
\end{array}, \begin{smallmatrix}
 0 &  1 &  0 &  0 \\
 1 &  0 &  0 &  0 \\
 0 &  0 &  0 &  1 \\
 0 &  0 &  1 &  0
\end{smallmatrix}\right)  & \leftrightarrow &
\left(\begin{smallmatrix}
 1 &  0 &  0 &  0 \\
 0 &  1 &  0 &  0 \\
 0 &  0 &  1 &  0 \\
 0 &  0 &  0 &  1
\end{smallmatrix}, \begin{array}{ccc}
4 & 2 & 2
\end{array}\right) \\[0.4cm]
 \left(\begin{array}{cc}
3 & 3
\end{array}, \begin{smallmatrix}
 0 &  1 &  0 &  0 \\
 1 & -1 &  1 &  0 \\
 0 &  1 & -1 &  1 \\
 0 &  0 &  1 &  0
\end{smallmatrix}\right)  & \leftrightarrow &
\left(\begin{smallmatrix}
 1 &  0 &  0 &  0 \\
 0 &  1 &  0 &  0 \\
 0 &  0 &  1 &  0 \\
 0 &  0 &  0 &  1
\end{smallmatrix}, \begin{array}{ccc}
4 & 3 & 2 \\
  & 2
\end{array}\right) &
\left(\begin{array}{cc}
3 & 3
\end{array}, \begin{smallmatrix}
 0 &  1 &  0 &  0 \\
 1 & -1 &  0 &  1 \\
 0 &  1 &  0 &  0 \\
 0 &  0 &  1 &  0
\end{smallmatrix}\right)  & \leftrightarrow &
\left(\begin{smallmatrix}
 1 &  0 &  0 &  0 \\
 0 &  1 &  0 &  0 \\
 0 &  0 &  1 &  0 \\
 0 &  0 &  0 &  1
\end{smallmatrix}, \begin{array}{ccc}
4 & 3 & 3 \\
  & 2
\end{array}\right) & \left(\begin{array}{cc}
3 & 3
\end{array}, \begin{smallmatrix}
 0 &  1 &  0 &  0 \\
 0 &  0 &  1 &  0 \\
 1 &  0 & -1 &  1 \\
 0 &  0 &  1 &  0
\end{smallmatrix}\right)  & \leftrightarrow &
\left(\begin{smallmatrix}
 1 &  0 &  0 &  0 \\
 0 &  0 &  1 &  0 \\
 0 &  0 &  0 &  1 \\
 0 &  1 &  0 &  0
\end{smallmatrix}, \begin{array}{ccc}
4 & 3 & 3 \\
  & 2
\end{array}\right) & \left(\begin{array}{cc}
3 & 3
\end{array}, \begin{smallmatrix}
 0 &  1 &  0 &  0 \\
 0 &  0 &  0 &  1 \\
 1 &  0 &  0 &  0 \\
 0 &  0 &  1 &  0
\end{smallmatrix}\right)  & \leftrightarrow &
\left(\begin{smallmatrix}
 1 &  0 &  0 &  0 \\
 0 &  0 &  1 &  0 \\
 0 &  1 & -1 &  1 \\
 0 &  0 &  1 &  0
\end{smallmatrix}, \begin{array}{ccc}
4 & 3 & 3 \\
  & 2
\end{array}\right)\\[0.4cm]
\left(\begin{array}{cc}
3 & 3
\end{array}, \begin{smallmatrix}
 0 &  1 &  0 &  0 \\
 1 & -1 &  1 &  0 \\
 0 &  0 &  0 &  1 \\
 0 &  1 &  0 &  0
\end{smallmatrix}\right)  & \leftrightarrow &
\left(\begin{smallmatrix}
 1 &  0 &  0 &  0 \\
 0 &  0 &  1 &  0 \\
 0 &  1 &  0 &  0 \\
 0 &  0 &  0 &  1
\end{smallmatrix}, \begin{array}{ccc}
4 & 3 & 2 \\
  & 2
\end{array}\right) & \left(\begin{array}{cc}
3 & 3
\end{array}, \begin{smallmatrix}
 0 &  1 &  0 &  0 \\
 1 & -1 &  0 &  1 \\
 0 &  0 &  1 &  0 \\
 0 &  1 &  0 &  0
\end{smallmatrix}\right)  & \leftrightarrow &
\left(\begin{smallmatrix}
 1 &  0 &  0 &  0 \\
 0 &  0 &  0 &  1 \\
 0 &  0 &  1 &  0 \\
 0 &  1 &  0 &  0
\end{smallmatrix}, \begin{array}{ccc}
4 & 3 & 3 \\
  & 2
\end{array}\right) & \left(\begin{array}{cc}
3 & 3
\end{array}, \begin{smallmatrix}
 0 &  1 &  0 &  0 \\
 0 &  0 &  1 &  0 \\
 1 & -1 &  0 &  1 \\
 0 &  1 &  0 &  0
\end{smallmatrix}\right)  & \leftrightarrow &
\left(\begin{smallmatrix}
 1 &  0 &  0 &  0 \\
 0 &  1 &  0 &  0 \\
 0 &  0 &  0 &  1 \\
 0 &  0 &  1 &  0
\end{smallmatrix}, \begin{array}{ccc}
4 & 3 & 1 \\
  & 2
\end{array}\right) &
\left(\begin{array}{cc}
3 & 3
\end{array}, \begin{smallmatrix}
 0 &  1 &  0 &  0 \\
 0 &  0 &  0 &  1 \\
 1 & -1 &  1 &  0 \\
 0 &  1 &  0 &  0
\end{smallmatrix}\right)  & \leftrightarrow &
\left(\begin{smallmatrix}
 1 &  0 &  0 &  0 \\
 0 &  1 &  0 &  0 \\
 0 &  0 &  0 &  1 \\
 0 &  0 &  1 &  0
\end{smallmatrix}, \begin{array}{ccc}
4 & 3 & 3 \\
  & 2
\end{array}\right) \\[0.4cm]
\left(\begin{array}{cc}
3 & 3
\end{array}, \begin{smallmatrix}
 0 &  1 &  0 &  0 \\
 0 &  0 &  1 &  0 \\
 0 &  0 &  0 &  1 \\
 1 &  0 &  0 &  0
\end{smallmatrix}\right)  & \leftrightarrow &
\left(\begin{smallmatrix}
 1 &  0 &  0 &  0 \\
 0 &  0 &  1 &  0 \\
 0 &  1 & -1 &  1 \\
 0 &  0 &  1 &  0
\end{smallmatrix}, \begin{array}{ccc}
4 & 3 & 2 \\
  & 2
\end{array}\right) & \left(\begin{array}{cc}
3 & 3
\end{array}, \begin{smallmatrix}
 0 &  1 &  0 &  0 \\
 0 &  0 &  0 &  1 \\
 0 &  0 &  1 &  0 \\
 1 &  0 &  0 &  0
\end{smallmatrix}\right)  & \leftrightarrow &
\left(\begin{smallmatrix}
 1 &  0 &  0 &  0 \\
 0 &  1 &  0 &  0 \\
 0 &  0 &  0 &  1 \\
 0 &  0 &  1 &  0
\end{smallmatrix}, \begin{array}{ccc}
4 & 3 & 2 \\
  & 2
\end{array}\right) &
\left(\begin{array}{cc}
3 & 3 \\
  & 2
\end{array}, \begin{smallmatrix}
 0 &  1 &  0 &  0 \\
 1 &  0 &  0 &  0 \\
 0 &  0 &  1 &  0 \\
 0 &  0 &  0 &  1
\end{smallmatrix}\right)  & \leftrightarrow &
\left(\begin{smallmatrix}
 1 &  0 &  0 &  0 \\
 0 &  1 &  0 &  0 \\
 0 &  0 &  1 &  0 \\
 0 &  0 &  0 &  1
\end{smallmatrix}, \begin{array}{cc}
4 & 3 \\
  & 2
\end{array}\right) & \left(\begin{array}{cc}
3 & 3 \\
  & 2
\end{array}, \begin{smallmatrix}
 0 &  1 &  0 &  0 \\
 1 & -1 &  1 &  0 \\
 0 &  1 &  0 &  0 \\
 0 &  0 &  0 &  1
\end{smallmatrix}\right)  & \leftrightarrow &
\left(\begin{smallmatrix}
 1 &  0 &  0 &  0 \\
 0 &  0 &  1 &  0 \\
 0 &  1 &  0 &  0 \\
 0 &  0 &  0 &  1
\end{smallmatrix}, \begin{array}{cc}
4 & 3 \\
  & 2
\end{array}\right) \\[0.4cm]
\left(\begin{array}{cc}
3 & 3 \\
  & 2
\end{array}, \begin{smallmatrix}
 0 &  1 &  0 &  0 \\
 0 &  0 &  1 &  0 \\
 1 &  0 &  0 &  0 \\
 0 &  0 &  0 &  1
\end{smallmatrix}\right)  & \leftrightarrow &
\left(\begin{smallmatrix}
 1 &  0 &  0 &  0 \\
 0 &  0 &  1 &  0 \\
 0 &  1 &  0 &  0 \\
 0 &  0 &  0 &  1
\end{smallmatrix}, \begin{array}{ccc}
4 & 3 & 1 \\
  & 2
\end{array}\right) &
\left(\begin{array}{cc}
3 & 3 \\
  & 2
\end{array}, \begin{smallmatrix}
 0 &  1 &  0 &  0 \\
 1 &  0 &  0 &  0 \\
 0 &  0 &  0 &  1 \\
 0 &  0 &  1 &  0
\end{smallmatrix}\right)  & \leftrightarrow &
\left(\begin{smallmatrix}
 1 &  0 &  0 &  0 \\
 0 &  1 &  0 &  0 \\
 0 &  0 &  0 &  1 \\
 0 &  0 &  1 &  0
\end{smallmatrix}, \begin{array}{cc}
4 & 3 \\
  & 2
\end{array}\right) & \left(\begin{array}{cc}
3 & 3 \\
  & 2
\end{array}, \begin{smallmatrix}
 0 &  1 &  0 &  0 \\
 1 & -1 &  1 &  0 \\
 0 &  1 & -1 &  1 \\
 0 &  0 &  1 &  0
\end{smallmatrix}\right)  & \leftrightarrow &
\left(\begin{smallmatrix}
 1 &  0 &  0 &  0 \\
 0 &  0 &  1 &  0 \\
 0 &  1 & -1 &  1 \\
 0 &  0 &  1 &  0
\end{smallmatrix}, \begin{array}{cc}
4 & 3 \\
  & 2
\end{array}\right) & \left(\begin{array}{cc}
3 & 3 \\
  & 2
\end{array}, \begin{smallmatrix}
 0 &  1 &  0 &  0 \\
 1 & -1 &  0 &  1 \\
 0 &  1 &  0 &  0 \\
 0 &  0 &  1 &  0
\end{smallmatrix}\right)  & \leftrightarrow &
\left(\begin{smallmatrix}
 1 &  0 &  0 &  0 \\
 0 &  0 &  0 &  1 \\
 0 &  1 &  0 &  0 \\
 0 &  0 &  1 &  0
\end{smallmatrix}, \begin{array}{cc}
4 & 3 \\
  & 2
\end{array}\right)\\[0.4cm]
\left(\begin{array}{cc}
3 & 3 \\
  & 2
\end{array}, \begin{smallmatrix}
 0 &  1 &  0 &  0 \\
 0 &  0 &  1 &  0 \\
 1 &  0 & -1 &  1 \\
 0 &  0 &  1 &  0
\end{smallmatrix}\right)  & \leftrightarrow &
\left(\begin{smallmatrix}
 1 &  0 &  0 &  0 \\
 0 &  0 &  1 &  0 \\
 0 &  1 & -1 &  1 \\
 0 &  0 &  1 &  0
\end{smallmatrix}, \begin{array}{ccc}
4 & 3 & 1 \\
  & 2
\end{array}\right) & \left(\begin{array}{cc}
3 & 3 \\
  & 2
\end{array}, \begin{smallmatrix}
 0 &  1 &  0 &  0 \\
 0 &  0 &  0 &  1 \\
 1 &  0 &  0 &  0 \\
 0 &  0 &  1 &  0
\end{smallmatrix}\right)  & \leftrightarrow &
\left(\begin{smallmatrix}
 1 &  0 &  0 &  0 \\
 0 &  0 &  0 &  1 \\
 0 &  1 &  0 &  0 \\
 0 &  0 &  1 &  0
\end{smallmatrix}, \begin{array}{ccc}
4 & 3 & 1 \\
  & 2
\end{array}\right) & \left(\begin{array}{cc}
3 & 3 \\
  & 2
\end{array}, \begin{smallmatrix}
 0 &  1 &  0 &  0 \\
 1 & -1 &  1 &  0 \\
 0 &  0 &  0 &  1 \\
 0 &  1 &  0 &  0
\end{smallmatrix}\right)  & \leftrightarrow &
\left(\begin{smallmatrix}
 1 &  0 &  0 &  0 \\
 0 &  0 &  1 &  0 \\
 0 &  0 &  0 &  1 \\
 0 &  1 &  0 &  0
\end{smallmatrix}, \begin{array}{cc}
4 & 3 \\
  & 2
\end{array}\right) &
\left(\begin{array}{cc}
3 & 3 \\
  & 2
\end{array}, \begin{smallmatrix}
 0 &  1 &  0 &  0 \\
 1 & -1 &  0 &  1 \\
 0 &  0 &  1 &  0 \\
 0 &  1 &  0 &  0
\end{smallmatrix}\right)  & \leftrightarrow &
\left(\begin{smallmatrix}
 1 &  0 &  0 &  0 \\
 0 &  0 &  0 &  1 \\
 0 &  0 &  1 &  0 \\
 0 &  1 &  0 &  0
\end{smallmatrix}, \begin{array}{cc}
4 & 3 \\
  & 2
\end{array}\right) \\[0.4cm]
\left(\begin{array}{cc}
3 & 3 \\
  & 2
\end{array}, \begin{smallmatrix}
 0 &  1 &  0 &  0 \\
 0 &  0 &  1 &  0 \\
 1 & -1 &  0 &  1 \\
 0 &  1 &  0 &  0
\end{smallmatrix}\right)  & \leftrightarrow &
\left(\begin{smallmatrix}
 1 &  0 &  0 &  0 \\
 0 &  0 &  1 &  0 \\
 0 &  0 &  0 &  1 \\
 0 &  1 &  0 &  0
\end{smallmatrix}, \begin{array}{ccc}
4 & 3 & 1 \\
  & 2
\end{array}\right) & \left(\begin{array}{cc}
3 & 3 \\
  & 2
\end{array}, \begin{smallmatrix}
 0 &  1 &  0 &  0 \\
 0 &  0 &  0 &  1 \\
 1 & -1 &  1 &  0 \\
 0 &  1 &  0 &  0
\end{smallmatrix}\right)  & \leftrightarrow &
\left(\begin{smallmatrix}
 1 &  0 &  0 &  0 \\
 0 &  0 &  0 &  1 \\
 0 &  0 &  1 &  0 \\
 0 &  1 &  0 &  0
\end{smallmatrix}, \begin{array}{ccc}
4 & 3 & 1 \\
  & 2
\end{array}\right) &
\left(\begin{array}{cc}
3 & 3 \\
  & 2
\end{array}, \begin{smallmatrix}
 0 &  1 &  0 &  0 \\
 0 &  0 &  1 &  0 \\
 0 &  0 &  0 &  1 \\
 1 &  0 &  0 &  0
\end{smallmatrix}\right)  & \leftrightarrow &
\left(\begin{smallmatrix}
 1 &  0 &  0 &  0 \\
 0 &  0 &  1 &  0 \\
 0 &  0 &  0 &  1 \\
 0 &  1 &  0 &  0
\end{smallmatrix}, \begin{array}{ccc}
4 & 3 & 2 \\
  & 2
\end{array}\right) & \left(\begin{array}{cc}
3 & 3 \\
  & 2
\end{array}, \begin{smallmatrix}
 0 &  1 &  0 &  0 \\
 0 &  0 &  0 &  1 \\
 0 &  0 &  1 &  0 \\
 1 &  0 &  0 &  0
\end{smallmatrix}\right)  & \leftrightarrow &
\left(\begin{smallmatrix}
 1 &  0 &  0 &  0 \\
 0 &  0 &  0 &  1 \\
 0 &  0 &  1 &  0 \\
 0 &  1 &  0 &  0
\end{smallmatrix}, \begin{array}{ccc}
4 & 3 & 2 \\
  & 2
\end{array}\right) \\[0.4cm]
\left(\begin{array}{cc}
3 & 2
\end{array}, \begin{smallmatrix}
 0 &  1 &  0 &  0 \\
 1 &  0 &  0 &  0 \\
 0 &  0 &  1 &  0 \\
 0 &  0 &  0 &  1
\end{smallmatrix}\right)  & \leftrightarrow &
\left(\begin{smallmatrix}
 1 &  0 &  0 &  0 \\
 0 &  1 &  0 &  0 \\
 0 &  0 &  1 &  0 \\
 0 &  0 &  0 &  1
\end{smallmatrix}, \begin{array}{ccc}
4 & 3 & 3
\end{array}\right) &
\left(\begin{array}{cc}
3 & 2
\end{array}, \begin{smallmatrix}
 0 &  1 &  0 &  0 \\
 1 & -1 &  1 &  0 \\
 0 &  1 &  0 &  0 \\
 0 &  0 &  0 &  1
\end{smallmatrix}\right)  & \leftrightarrow &
\left(\begin{smallmatrix}
 1 &  0 &  0 &  0 \\
 0 &  0 &  1 &  0 \\
 0 &  1 &  0 &  0 \\
 0 &  0 &  0 &  1
\end{smallmatrix}, \begin{array}{ccc}
4 & 3 & 3
\end{array}\right) & \left(\begin{array}{cc}
3 & 2
\end{array}, \begin{smallmatrix}
 0 &  1 &  0 &  0 \\
 0 &  0 &  1 &  0 \\
 1 &  0 &  0 &  0 \\
 0 &  0 &  0 &  1
\end{smallmatrix}\right)  & \leftrightarrow &
\left(\begin{smallmatrix}
 1 &  0 &  0 &  0 \\
 0 &  0 &  1 &  0 \\
 0 &  1 & -1 &  1 \\
 0 &  0 &  1 &  0
\end{smallmatrix}, \begin{array}{c}
4
\end{array}\right) & \left(\begin{array}{cc}
3 & 2
\end{array}, \begin{smallmatrix}
 0 &  1 &  0 &  0 \\
 1 &  0 &  0 &  0 \\
 0 &  0 &  0 &  1 \\
 0 &  0 &  1 &  0
\end{smallmatrix}\right)  & \leftrightarrow &
\left(\begin{smallmatrix}
 1 &  0 &  0 &  0 \\
 0 &  1 &  0 &  0 \\
 0 &  0 &  0 &  1 \\
 0 &  0 &  1 &  0
\end{smallmatrix}, \begin{array}{ccc}
4 & 3 & 3
\end{array}\right)\\[0.4cm]
\left(\begin{array}{cc}
3 & 2
\end{array}, \begin{smallmatrix}
 0 &  1 &  0 &  0 \\
 1 & -1 &  1 &  0 \\
 0 &  1 & -1 &  1 \\
 0 &  0 &  1 &  0
\end{smallmatrix}\right)  & \leftrightarrow &
\left(\begin{smallmatrix}
 1 &  0 &  0 &  0 \\
 0 &  0 &  1 &  0 \\
 0 &  1 & -1 &  1 \\
 0 &  0 &  1 &  0
\end{smallmatrix}, \begin{array}{ccc}
4 & 3 & 3
\end{array}\right) & \left(\begin{array}{cc}
3 & 2
\end{array}, \begin{smallmatrix}
 0 &  1 &  0 &  0 \\
 1 & -1 &  0 &  1 \\
 0 &  1 &  0 &  0 \\
 0 &  0 &  1 &  0
\end{smallmatrix}\right)  & \leftrightarrow &
\left(\begin{smallmatrix}
 1 &  0 &  0 &  0 \\
 0 &  0 &  0 &  1 \\
 0 &  1 &  0 &  0 \\
 0 &  0 &  1 &  0
\end{smallmatrix}, \begin{array}{ccc}
4 & 3 & 3
\end{array}\right) & \left(\begin{array}{cc}
3 & 2
\end{array}, \begin{smallmatrix}
 0 &  1 &  0 &  0 \\
 0 &  0 &  1 &  0 \\
 1 &  0 & -1 &  1 \\
 0 &  0 &  1 &  0
\end{smallmatrix}\right)  & \leftrightarrow &
\left(\begin{smallmatrix}
 1 &  0 &  0 &  0 \\
 0 &  1 &  0 &  0 \\
 0 &  0 &  0 &  1 \\
 0 &  0 &  1 &  0
\end{smallmatrix}, \begin{array}{ccc}
4 & 2 & 2
\end{array}\right) &
\left(\begin{array}{cc}
3 & 2
\end{array}, \begin{smallmatrix}
 0 &  1 &  0 &  0 \\
 0 &  0 &  0 &  1 \\
 1 &  0 &  0 &  0 \\
 0 &  0 &  1 &  0
\end{smallmatrix}\right)  & \leftrightarrow &
\left(\begin{smallmatrix}
 1 &  0 &  0 &  0 \\
 0 &  0 &  0 &  1 \\
 0 &  1 &  0 &  0 \\
 0 &  0 &  1 &  0
\end{smallmatrix}, \begin{array}{ccc}
4 & 1 & 1
\end{array}\right) \\[0.4cm]
 \left(\begin{array}{cc}
3 & 2
\end{array}, \begin{smallmatrix}
 0 &  1 &  0 &  0 \\
 1 & -1 &  1 &  0 \\
 0 &  0 &  0 &  1 \\
 0 &  1 &  0 &  0
\end{smallmatrix}\right)  & \leftrightarrow &
\left(\begin{smallmatrix}
 1 &  0 &  0 &  0 \\
 0 &  0 &  1 &  0 \\
 0 &  0 &  0 &  1 \\
 0 &  1 &  0 &  0
\end{smallmatrix}, \begin{array}{ccc}
4 & 3 & 3
\end{array}\right) & \left(\begin{array}{cc}
3 & 2
\end{array}, \begin{smallmatrix}
 0 &  1 &  0 &  0 \\
 1 & -1 &  0 &  1 \\
 0 &  0 &  1 &  0 \\
 0 &  1 &  0 &  0
\end{smallmatrix}\right)  & \leftrightarrow &
\left(\begin{smallmatrix}
 1 &  0 &  0 &  0 \\
 0 &  0 &  0 &  1 \\
 0 &  0 &  1 &  0 \\
 0 &  1 &  0 &  0
\end{smallmatrix}, \begin{array}{ccc}
4 & 3 & 3
\end{array}\right) &
\left(\begin{array}{cc}
3 & 2
\end{array}, \begin{smallmatrix}
 0 &  1 &  0 &  0 \\
 0 &  0 &  1 &  0 \\
 1 & -1 &  0 &  1 \\
 0 &  1 &  0 &  0
\end{smallmatrix}\right)  & \leftrightarrow &
\left(\begin{smallmatrix}
 1 &  0 &  0 &  0 \\
 0 &  0 &  1 &  0 \\
 0 &  1 &  0 &  0 \\
 0 &  0 &  0 &  1
\end{smallmatrix}, \begin{array}{ccc}
4 & 1 & 1
\end{array}\right) & \left(\begin{array}{cc}
3 & 2
\end{array}, \begin{smallmatrix}
 0 &  1 &  0 &  0 \\
 0 &  0 &  0 &  1 \\
 1 & -1 &  1 &  0 \\
 0 &  1 &  0 &  0
\end{smallmatrix}\right)  & \leftrightarrow &
\left(\begin{smallmatrix}
 1 &  0 &  0 &  0 \\
 0 &  0 &  1 &  0 \\
 0 &  1 &  0 &  0 \\
 0 &  0 &  0 &  1
\end{smallmatrix}, \begin{array}{ccc}
4 & 3 & 3 \\
  & 2
\end{array}\right) \\[0.4cm]
\left(\begin{array}{cc}
3 & 2
\end{array}, \begin{smallmatrix}
 0 &  1 &  0 &  0 \\
 0 &  0 &  1 &  0 \\
 0 &  0 &  0 &  1 \\
 1 &  0 &  0 &  0
\end{smallmatrix}\right)  & \leftrightarrow &
\left(\begin{smallmatrix}
 1 &  0 &  0 &  0 \\
 0 &  0 &  1 &  0 \\
 0 &  1 &  0 &  0 \\
 0 &  0 &  0 &  1
\end{smallmatrix}, \begin{array}{cc}
4 & 1
\end{array}\right) &
\left(\begin{array}{cc}
3 & 2
\end{array}, \begin{smallmatrix}
 0 &  1 &  0 &  0 \\
 0 &  0 &  0 &  1 \\
 0 &  0 &  1 &  0 \\
 1 &  0 &  0 &  0
\end{smallmatrix}\right)  & \leftrightarrow &
\left(\begin{smallmatrix}
 1 &  0 &  0 &  0 \\
 0 &  0 &  0 &  1 \\
 0 &  1 &  0 &  0 \\
 0 &  0 &  1 &  0
\end{smallmatrix}, \begin{array}{ccc}
4 & 3 & 3 \\
  & 2
\end{array}\right) & \left(\begin{array}{cc}
3 & 1
\end{array}, \begin{smallmatrix}
 0 &  1 &  0 &  0 \\
 1 &  0 &  0 &  0 \\
 0 &  0 &  1 &  0 \\
 0 &  0 &  0 &  1
\end{smallmatrix}\right)  & \leftrightarrow &
\left(\begin{smallmatrix}
 1 &  0 &  0 &  0 \\
 0 &  1 &  0 &  0 \\
 0 &  0 &  1 &  0 \\
 0 &  0 &  0 &  1
\end{smallmatrix}, \begin{array}{ccc}
4 & 3 & 2
\end{array}\right) & \left(\begin{array}{cc}
3 & 1
\end{array}, \begin{smallmatrix}
 0 &  1 &  0 &  0 \\
 1 & -1 &  1 &  0 \\
 0 &  1 &  0 &  0 \\
 0 &  0 &  0 &  1
\end{smallmatrix}\right)  & \leftrightarrow &
\left(\begin{smallmatrix}
 1 &  0 &  0 &  0 \\
 0 &  0 &  1 &  0 \\
 0 &  1 &  0 &  0 \\
 0 &  0 &  0 &  1
\end{smallmatrix}, \begin{array}{ccc}
4 & 3 & 2
\end{array}\right)\\[0.4cm]
\left(\begin{array}{cc}
3 & 1
\end{array}, \begin{smallmatrix}
 0 &  1 &  0 &  0 \\
 0 &  0 &  1 &  0 \\
 1 &  0 &  0 &  0 \\
 0 &  0 &  0 &  1
\end{smallmatrix}\right)  & \leftrightarrow &
\left(\begin{smallmatrix}
 1 &  0 &  0 &  0 \\
 0 &  0 &  1 &  0 \\
 0 &  1 &  0 &  0 \\
 0 &  0 &  0 &  1
\end{smallmatrix}, \begin{array}{ccc}
4 & 2 & 2
\end{array}\right) & \left(\begin{array}{cc}
3 & 1
\end{array}, \begin{smallmatrix}
 0 &  1 &  0 &  0 \\
 1 &  0 &  0 &  0 \\
 0 &  0 &  0 &  1 \\
 0 &  0 &  1 &  0
\end{smallmatrix}\right)  & \leftrightarrow &
\left(\begin{smallmatrix}
 1 &  0 &  0 &  0 \\
 0 &  1 &  0 &  0 \\
 0 &  0 &  0 &  1 \\
 0 &  0 &  1 &  0
\end{smallmatrix}, \begin{array}{ccc}
4 & 3 & 2
\end{array}\right) & \left(\begin{array}{cc}
3 & 1
\end{array}, \begin{smallmatrix}
 0 &  1 &  0 &  0 \\
 1 & -1 &  1 &  0 \\
 0 &  1 & -1 &  1 \\
 0 &  0 &  1 &  0
\end{smallmatrix}\right)  & \leftrightarrow &
\left(\begin{smallmatrix}
 1 &  0 &  0 &  0 \\
 0 &  0 &  1 &  0 \\
 0 &  1 & -1 &  1 \\
 0 &  0 &  1 &  0
\end{smallmatrix}, \begin{array}{ccc}
4 & 3 & 2
\end{array}\right) &
\left(\begin{array}{cc}
3 & 1
\end{array}, \begin{smallmatrix}
 0 &  1 &  0 &  0 \\
 1 & -1 &  0 &  1 \\
 0 &  1 &  0 &  0 \\
 0 &  0 &  1 &  0
\end{smallmatrix}\right)  & \leftrightarrow &
\left(\begin{smallmatrix}
 1 &  0 &  0 &  0 \\
 0 &  0 &  0 &  1 \\
 0 &  1 &  0 &  0 \\
 0 &  0 &  1 &  0
\end{smallmatrix}, \begin{array}{ccc}
4 & 3 & 2
\end{array}\right) \\[0.4cm]
\left(\begin{array}{cc}
3 & 1
\end{array}, \begin{smallmatrix}
 0 &  1 &  0 &  0 \\
 0 &  0 &  1 &  0 \\
 1 &  0 & -1 &  1 \\
 0 &  0 &  1 &  0
\end{smallmatrix}\right)  & \leftrightarrow &
\left(\begin{smallmatrix}
 1 &  0 &  0 &  0 \\
 0 &  0 &  1 &  0 \\
 0 &  1 & -1 &  1 \\
 0 &  0 &  1 &  0
\end{smallmatrix}, \begin{array}{ccc}
4 & 2 & 2
\end{array}\right) & \left(\begin{array}{cc}
3 & 1
\end{array}, \begin{smallmatrix}
 0 &  1 &  0 &  0 \\
 0 &  0 &  0 &  1 \\
 1 &  0 &  0 &  0 \\
 0 &  0 &  1 &  0
\end{smallmatrix}\right)  & \leftrightarrow &
\left(\begin{smallmatrix}
 1 &  0 &  0 &  0 \\
 0 &  0 &  0 &  1 \\
 0 &  1 &  0 &  0 \\
 0 &  0 &  1 &  0
\end{smallmatrix}, \begin{array}{ccc}
4 & 2 & 2
\end{array}\right) &
\left(\begin{array}{cc}
3 & 1
\end{array}, \begin{smallmatrix}
 0 &  1 &  0 &  0 \\
 1 & -1 &  1 &  0 \\
 0 &  0 &  0 &  1 \\
 0 &  1 &  0 &  0
\end{smallmatrix}\right)  & \leftrightarrow &
\left(\begin{smallmatrix}
 1 &  0 &  0 &  0 \\
 0 &  0 &  1 &  0 \\
 0 &  0 &  0 &  1 \\
 0 &  1 &  0 &  0
\end{smallmatrix}, \begin{array}{ccc}
4 & 3 & 2
\end{array}\right) & \left(\begin{array}{cc}
3 & 1
\end{array}, \begin{smallmatrix}
 0 &  1 &  0 &  0 \\
 1 & -1 &  0 &  1 \\
 0 &  0 &  1 &  0 \\
 0 &  1 &  0 &  0
\end{smallmatrix}\right)  & \leftrightarrow &
\left(\begin{smallmatrix}
 1 &  0 &  0 &  0 \\
 0 &  0 &  0 &  1 \\
 0 &  0 &  1 &  0 \\
 0 &  1 &  0 &  0
\end{smallmatrix}, \begin{array}{ccc}
4 & 3 & 2
\end{array}\right) \\[0.4cm]
\left(\begin{array}{cc}
3 & 1
\end{array}, \begin{smallmatrix}
 0 &  1 &  0 &  0 \\
 0 &  0 &  1 &  0 \\
 1 & -1 &  0 &  1 \\
 0 &  1 &  0 &  0
\end{smallmatrix}\right)  & \leftrightarrow &
\left(\begin{smallmatrix}
 1 &  0 &  0 &  0 \\
 0 &  0 &  1 &  0 \\
 0 &  0 &  0 &  1 \\
 0 &  1 &  0 &  0
\end{smallmatrix}, \begin{array}{ccc}
4 & 2 & 2
\end{array}\right) &
\left(\begin{array}{cc}
3 & 1
\end{array}, \begin{smallmatrix}
 0 &  1 &  0 &  0 \\
 0 &  0 &  0 &  1 \\
 1 & -1 &  1 &  0 \\
 0 &  1 &  0 &  0
\end{smallmatrix}\right)  & \leftrightarrow &
\left(\begin{smallmatrix}
 1 &  0 &  0 &  0 \\
 0 &  0 &  0 &  1 \\
 0 &  0 &  1 &  0 \\
 0 &  1 &  0 &  0
\end{smallmatrix}, \begin{array}{ccc}
4 & 2 & 2
\end{array}\right) & \left(\begin{array}{cc}
3 & 1
\end{array}, \begin{smallmatrix}
 0 &  1 &  0 &  0 \\
 0 &  0 &  1 &  0 \\
 0 &  0 &  0 &  1 \\
 1 &  0 &  0 &  0
\end{smallmatrix}\right)  & \leftrightarrow &
\left(\begin{smallmatrix}
 1 &  0 &  0 &  0 \\
 0 &  0 &  1 &  0 \\
 0 &  1 &  0 &  0 \\
 0 &  0 &  0 &  1
\end{smallmatrix}, \begin{array}{c}
4
\end{array}\right) & \left(\begin{array}{cc}
3 & 1
\end{array}, \begin{smallmatrix}
 0 &  1 &  0 &  0 \\
 0 &  0 &  0 &  1 \\
 0 &  0 &  1 &  0 \\
 1 &  0 &  0 &  0
\end{smallmatrix}\right)  & \leftrightarrow &
\left(\begin{smallmatrix}
 1 &  0 &  0 &  0 \\
 0 &  0 &  0 &  1 \\
 0 &  1 &  0 &  0 \\
 0 &  0 &  1 &  0
\end{smallmatrix}, \begin{array}{c}
4
\end{array}\right)\\[0.4cm]
\left(\begin{array}{c}
3
\end{array}, \begin{smallmatrix}
 0 &  1 &  0 &  0 \\
 1 &  0 &  0 &  0 \\
 0 &  0 &  1 &  0 \\
 0 &  0 &  0 &  1
\end{smallmatrix}\right)  & \leftrightarrow &
\left(\begin{smallmatrix}
 1 &  0 &  0 &  0 \\
 0 &  1 &  0 &  0 \\
 0 &  0 &  1 &  0 \\
 0 &  0 &  0 &  1
\end{smallmatrix}, \begin{array}{ccc}
4 & 3 & 1
\end{array}\right) & \left(\begin{array}{c}
3
\end{array}, \begin{smallmatrix}
 0 &  1 &  0 &  0 \\
 1 & -1 &  1 &  0 \\
 0 &  1 &  0 &  0 \\
 0 &  0 &  0 &  1
\end{smallmatrix}\right)  & \leftrightarrow &
\left(\begin{smallmatrix}
 1 &  0 &  0 &  0 \\
 0 &  0 &  1 &  0 \\
 0 &  1 &  0 &  0 \\
 0 &  0 &  0 &  1
\end{smallmatrix}, \begin{array}{ccc}
4 & 3 & 1
\end{array}\right) & \left(\begin{array}{c}
3
\end{array}, \begin{smallmatrix}
 0 &  1 &  0 &  0 \\
 0 &  0 &  1 &  0 \\
 1 &  0 &  0 &  0 \\
 0 &  0 &  0 &  1
\end{smallmatrix}\right)  & \leftrightarrow &
\left(\begin{smallmatrix}
 1 &  0 &  0 &  0 \\
 0 &  0 &  1 &  0 \\
 0 &  1 &  0 &  0 \\
 0 &  0 &  0 &  1
\end{smallmatrix}, \begin{array}{ccc}
4 & 2 & 1
\end{array}\right) &
\left(\begin{array}{c}
3
\end{array}, \begin{smallmatrix}
 0 &  1 &  0 &  0 \\
 1 &  0 &  0 &  0 \\
 0 &  0 &  0 &  1 \\
 0 &  0 &  1 &  0
\end{smallmatrix}\right)  & \leftrightarrow &
\left(\begin{smallmatrix}
 1 &  0 &  0 &  0 \\
 0 &  1 &  0 &  0 \\
 0 &  0 &  0 &  1 \\
 0 &  0 &  1 &  0
\end{smallmatrix}, \begin{array}{ccc}
4 & 3 & 1
\end{array}\right) \\[0.4cm]
\left(\begin{array}{c}
3
\end{array}, \begin{smallmatrix}
 0 &  1 &  0 &  0 \\
 1 & -1 &  1 &  0 \\
 0 &  1 & -1 &  1 \\
 0 &  0 &  1 &  0
\end{smallmatrix}\right)  & \leftrightarrow &
\left(\begin{smallmatrix}
 1 &  0 &  0 &  0 \\
 0 &  0 &  1 &  0 \\
 0 &  1 & -1 &  1 \\
 0 &  0 &  1 &  0
\end{smallmatrix}, \begin{array}{ccc}
4 & 3 & 1
\end{array}\right) & \left(\begin{array}{c}
3
\end{array}, \begin{smallmatrix}
 0 &  1 &  0 &  0 \\
 1 & -1 &  0 &  1 \\
 0 &  1 &  0 &  0 \\
 0 &  0 &  1 &  0
\end{smallmatrix}\right)  & \leftrightarrow &
\left(\begin{smallmatrix}
 1 &  0 &  0 &  0 \\
 0 &  0 &  0 &  1 \\
 0 &  1 &  0 &  0 \\
 0 &  0 &  1 &  0
\end{smallmatrix}, \begin{array}{ccc}
4 & 3 & 1
\end{array}\right) &
\left(\begin{array}{c}
3
\end{array}, \begin{smallmatrix}
 0 &  1 &  0 &  0 \\
 0 &  0 &  1 &  0 \\
 1 &  0 & -1 &  1 \\
 0 &  0 &  1 &  0
\end{smallmatrix}\right)  & \leftrightarrow &
\left(\begin{smallmatrix}
 1 &  0 &  0 &  0 \\
 0 &  0 &  1 &  0 \\
 0 &  1 & -1 &  1 \\
 0 &  0 &  1 &  0
\end{smallmatrix}, \begin{array}{ccc}
4 & 2 & 1
\end{array}\right) & \left(\begin{array}{c}
3
\end{array}, \begin{smallmatrix}
 0 &  1 &  0 &  0 \\
 0 &  0 &  0 &  1 \\
 1 &  0 &  0 &  0 \\
 0 &  0 &  1 &  0
\end{smallmatrix}\right)  & \leftrightarrow &
\left(\begin{smallmatrix}
 1 &  0 &  0 &  0 \\
 0 &  0 &  0 &  1 \\
 0 &  1 &  0 &  0 \\
 0 &  0 &  1 &  0
\end{smallmatrix}, \begin{array}{ccc}
4 & 2 & 1
\end{array}\right) \\[0.4cm]
 \left(\begin{array}{c}
3
\end{array}, \begin{smallmatrix}
 0 &  1 &  0 &  0 \\
 1 & -1 &  1 &  0 \\
 0 &  0 &  0 &  1 \\
 0 &  1 &  0 &  0
\end{smallmatrix}\right)  & \leftrightarrow &
\left(\begin{smallmatrix}
 1 &  0 &  0 &  0 \\
 0 &  0 &  1 &  0 \\
 0 &  0 &  0 &  1 \\
 0 &  1 &  0 &  0
\end{smallmatrix}, \begin{array}{ccc}
4 & 3 & 1
\end{array}\right) &
\left(\begin{array}{c}
3
\end{array}, \begin{smallmatrix}
 0 &  1 &  0 &  0 \\
 1 & -1 &  0 &  1 \\
 0 &  0 &  1 &  0 \\
 0 &  1 &  0 &  0
\end{smallmatrix}\right)  & \leftrightarrow &
\left(\begin{smallmatrix}
 1 &  0 &  0 &  0 \\
 0 &  0 &  0 &  1 \\
 0 &  0 &  1 &  0 \\
 0 &  1 &  0 &  0
\end{smallmatrix}, \begin{array}{ccc}
4 & 3 & 1
\end{array}\right) & \left(\begin{array}{c}
3
\end{array}, \begin{smallmatrix}
 0 &  1 &  0 &  0 \\
 0 &  0 &  1 &  0 \\
 1 & -1 &  0 &  1 \\
 0 &  1 &  0 &  0
\end{smallmatrix}\right)  & \leftrightarrow &
\left(\begin{smallmatrix}
 1 &  0 &  0 &  0 \\
 0 &  0 &  1 &  0 \\
 0 &  0 &  0 &  1 \\
 0 &  1 &  0 &  0
\end{smallmatrix}, \begin{array}{ccc}
4 & 2 & 1
\end{array}\right) & \left(\begin{array}{c}
3
\end{array}, \begin{smallmatrix}
 0 &  1 &  0 &  0 \\
 0 &  0 &  0 &  1 \\
 1 & -1 &  1 &  0 \\
 0 &  1 &  0 &  0
\end{smallmatrix}\right)  & \leftrightarrow &
\left(\begin{smallmatrix}
 1 &  0 &  0 &  0 \\
 0 &  0 &  0 &  1 \\
 0 &  0 &  1 &  0 \\
 0 &  1 &  0 &  0
\end{smallmatrix}, \begin{array}{ccc}
4 & 2 & 1
\end{array}\right)\\[0.4cm]
\left(\begin{array}{c}
3
\end{array}, \begin{smallmatrix}
 0 &  1 &  0 &  0 \\
 0 &  0 &  1 &  0 \\
 0 &  0 &  0 &  1 \\
 1 &  0 &  0 &  0
\end{smallmatrix}\right)  & \leftrightarrow &
\left(\begin{smallmatrix}
 1 &  0 &  0 &  0 \\
 0 &  0 &  1 &  0 \\
 0 &  0 &  0 &  1 \\
 0 &  1 &  0 &  0
\end{smallmatrix}, \begin{array}{ccc}
4 & 1 & 1
\end{array}\right) & \left(\begin{array}{c}
3
\end{array}, \begin{smallmatrix}
 0 &  1 &  0 &  0 \\
 0 &  0 &  0 &  1 \\
 0 &  0 &  1 &  0 \\
 1 &  0 &  0 &  0
\end{smallmatrix}\right)  & \leftrightarrow &
\left(\begin{smallmatrix}
 1 &  0 &  0 &  0 \\
 0 &  0 &  0 &  1 \\
 0 &  0 &  1 &  0 \\
 0 &  1 &  0 &  0
\end{smallmatrix}, \begin{array}{ccc}
4 & 1 & 1
\end{array}\right)\\[0.4cm]
\end{array}$}
\end{center}

The process of finding our bijections started by translating known ``computational'' proofs (see \cite{Fis06,Fis07} for the original proof and \cite{Fis16} for the most concise version so far) into  bijective proofs. However, the combinatorial point of view led to several interesting modifications so that the original proofs are not always apparent.

\subsection*{Structure of the paper} We start by a short introduction on signed sets and sijections in Section~\ref{sect:signedsets}. Note that these are a special case of traced monoidal categories, see \cite{traced}. In Section~\ref{sect:B} we apply some of these notions to obtain a sijection on $\B_{n,i}$,  which serves also as a good warm-up. In Section~\ref{sect:LinAlg} we ``sijectify'' known results on determinants such as the multiplicativity and Cramer's rule. In Section~\ref{sect:pre}, we summarize some prerequisites needed from \cite{PartI} before we finish the construction in the final two sections.

\section{Signed sets and sijections} \label{sect:signedsets}

Our constructions rely heavily on the notion of signed sets. We use the terminology that was developed in \cite{PartI}. We repeat the key points here, and refer the reader to \cite[\S 2]{PartI} for all the details and examples.

\medskip

A \emph{signed set} is a pair of disjoint finite sets: $\u S = (S^+,S^-)$ with $S^+ \cap S^- = \emptyset$. Equivalently, a signed set is a finite set $S$ together with a sign function  $\sign \colon S \to \{1,-1\}$, but we will mostly avoid the use of the sign function. Signed sets are usually underlined throughout the paper with the following exception: an ordinary set $S$ always induces a signed set $\u S = (S,\emptyset)$, and in this case we identify $\u S$ with $S$. We summarize related notions.

\begin{itemize}
\item The \emph{size} of a signed set $\u S$ is $|\u S| = |S^+| - |S^-|$.
\item The \emph{opposite} signed set of $\u S$ is $- \u S = (S^-,S^+)$.
\item  The \emph{Cartesian product} of signed sets $\u S$ and $\u T$ is
$$\u S \times \u T = (S^+ \times T^+ \cup S^- \times T^-,S^+ \times T^- \cup S^- \times T^+).$$
\item The \emph{disjoint union} of signed sets $\u S$ and $\u T$ is
$$\u S \sqcup \u T = (\u S \times (\{0\},\emptyset)) \cup (\u T \times (\{1\},\emptyset)).$$
\item The \emph{disjoint union of a family of signed sets} $\u S_t$ indexed with a signed set $\u T$ is
$$\bigsqcup_{t \in \u T} \u S_t = \bigcup_{t \in \u T} (\u S_t \times \u{\{t\}}).$$
Here $\u{\{t\}}$ is $(\{t\},\emptyset)$ if $t \in T^+$ and $(\emptyset,\{t\})$ if $t \in T^-$.
\end{itemize}
Most of the usual properties of Cartesian products and disjoint unions of ordinary sets extend to signed sets.

An important type of signed sets are signed intervals: for $a,b \in \Z$, define
$$\u{[a,b]} = \begin{cases} ([a,b],\emptyset) & \mbox{if } a \leq b \\ (\emptyset,[b+1,a-1]) & \mbox{if } a > b \end{cases}.$$
Here $[a,b]$ stands for the usual interval in $\mathbb{Z}$. The signed sets that are of relevance in this paper are usually constructed from signed intervals using Cartesian products and disjoint unions.

The generalization of bijections to signed sets is played by ``signed bijections'', which we call \emph{sijections}. A sijection $\varphi$ from $\u S$ to $\u T$,
$$\varphi \colon \u S \Rightarrow \u T,$$
is an involution on the set $(S^+ \cup S^-) \sqcup (T^+ \cup T^-)$ with the property $\varphi(S^+ \sqcup T^-) = S^- \sqcup T^+$. It follows that also $\varphi(S^- \sqcup T^+) = S^+ \sqcup T^-$. A sijection can also be thought of as a collection of a sign-reversing involution on a subset of $\u S$, a sign-reversing involution on a subset of $\u T$, and a {sign-preserving} matching between the remaining elements of $\u S$ with the remaining elements of $\u T$.
The existence of a sijection $\varphi \colon \u S \Rightarrow \u T$ clearly implies $|\u S| = |S^+| - |S^-| = |T^+| - |T^-| = |\u T|$.

In Proposition~2 of \cite{PartI} it is explained how to construct the Cartesian product and the disjoint union of sijections, and also how to compose two sijections using a variant of the Garsia-Milne involution principle. These constructions are fundamental for most of the constructions in this paper. It follows that the existence of a sijection between $\u S$ and $\u T$ is an equivalence relation; it is denoted by ``$\approx$''.

The sijection that is underlying most of our constructions is the following.

\begin{problem}(\cite[Problem 1]{PartI}) \label{prop:alpha}
 Given $a,b,c \in \Z$, construct a sijection
 $$\alpha = \alpha_{a,b,c} \colon \si a c \Rightarrow \si a b \sqcup \si{b+1}c.$$
\end{problem}

\begin{proof}[Construction]
 For $a \leq b \leq c$ and $c < b < a$, there is nothing to prove. For, say, $a \leq c < b$, we have
 $$\si a b \sqcup \si{b+1}c = (\si a c \sqcup \si{c+1} b) \sqcup \si{b+1}c = \si a c \sqcup (\si{c+1} b \sqcup (-\si {c+1}{b}))$$
 and since there is a sijection $\si{c+1} b \sqcup (-\si {c+1}{b}) \Rightarrow \u \emptyset$, we get a sijection $\si a b \sqcup \si{b+1}c \Rightarrow \si a c$. The cases $b < a \leq c$, $b \leq c < a$, and $c < a \leq b$ are analogous. \comment{Note that in all cases, $\si a b$ and  $\si{b+1}c$ are either disjoint or one set is contained in the other.}
\end{proof}

It has been speculated by combinatorialists working on bijections between ASMs and classes of plane partitions such as DPPs or totally symmetric self-complementary plane partitions that these bijections will involve some sort of \emph{jeu de taquin}. We are wondering whether the application of $\alpha_{a,b,c}$ has anything to do with such a jeu de taquin move.

In \cite{PartI}, the notion of elementary signed sets and normal sijections turned out to be useful to simplify the description of several sijections.

\begin{definition}
 A signed set $\u A$ is \emph{elementary of dimension $n$ and depth $0$} if its elements are in $\Z^n$. A signed set $\u A$ is \emph{elementary of dimension $n$ and depth $d$}, $d \geq 1$, if it is of the form
 $$\bigsqcup_{t \in \u T} \u S_t,$$
 where $\u T$ is a signed set, and $\u S_t$ are all signed sets of dimension $n$ and depth at most $d-1$, with the depth of at least one of them equal to $d-1$. A signed set $\u A$ is \emph{elementary of dimension $n$} if it is an elementary signed set of dimension $n$ and depth $d$ for some $d \in \N$.\\
 The \emph{projection map} on an elementary set of dimension $n$ is the map
 $$\xi \colon \u A \to \Z^n$$
 defined as follows. If the depth of $\u A$ is $0$, then $\xi$ is simply the inclusion map. Once $\xi$ is defined on elementary signed sets of depth $< d$, and the depth of $\u A$ is $d$, then $\u A = \bigsqcup_{t \in \u T} \u S_t$, where $\xi$ is defined on all $\u S_t$. Then define $\xi(s,t) = \xi(s)$ for $(s,t) \in \u A$.\\
 A sijection $\psi \colon \u T \Rightarrow \u{\wt T}$ between elementary signed sets $\u T$ and $\u{\wt T}$ of the same dimension is \emph{normal} if $\xi(\psi(t)) = \xi(t)$ for all $t \in \u T \sqcup \u{\wt T}$.
\end{definition}

The sijection $\alpha$ from Problem~\ref{prop:alpha} is normal.

The main reason normal sijections are important is that they give a very natural special case for disjoint unions of sijections (see  \cite[Proposition~2 (3)]{PartI}). Suppose that $\u T$ and $\u{\wt T}$ are elementary signed sets of dimension $n$, and that $\psi \colon \u T \Rightarrow \u{\wt T}$ is a normal sijection. Furthermore, suppose that we have a signed set $\u S_{\q k}$ for every $\q k \in \Z^n$. Then we have a sijection
$$\bigsqcup_{t \in \u T} \u S_{\xi(t)} \Rightarrow \bigsqcup_{t \in \u{\wt T}} \u S_{\xi(t)}.$$
Indeed, \cite[Proposition~2 (3)]{PartI} gives us a sijection provided that we have a sijection $\varphi_t \colon \u S_{\xi(t)} \Rightarrow \u S_{\xi(\psi(t))}$ satisfying $\varphi_{\psi(t)} = \varphi_t^{-1}$ for every $t \in \u T \sqcup \u{\wt T}$. But since $\xi(\psi(t)) = \xi(t)$, we can take $\varphi_t$ to be the identity.

\section{A sijection on $\B_{n,i}$} \label{sect:B}

For a set $A$ and $n \in \N$, denote by $\binom A n$ the set of all $n$-element subsets of $A$.%, by $\p P_A$ the set of all subsets of $A$, and by $\u{\p P}_A$ the signed set of all subsets of $A$, where a subset is positive if and only if it has an even number of elements.

\begin{problem} \label{prob:trinomial}
 For $a,b,c \in \N$, construct a bijection
 $$\binom{[a+b+c]}a \times \binom{[b+c]}b \longrightarrow \binom{[a+b+c]}b \times \binom{[a+c]}c.$$
\end{problem}
\begin{proof}[Construction]
 An element on the left can be interpreted as a pair $(A,B)$ of disjoint subsets of $[a+b+c]$ with $|A|=a$ and $|B|=b$, and an element on the right can be interpreted as a pair $(B,C)$ of disjoint subsets of $[a+b+c]$ with $|B|=b$ and $|C|=c$. The bijection is then $(A,B) \mapsto (B, [a+b+c] \setminus (A \cup B))$.
\end{proof}

The following construction gives a sijective proof of a variation of the Chu-Vandermonde identity.

\begin{problem} \label{prob:binom}
 For $a,b,c \in \N$, construct a sijection
 $$\bigsqcup_{j=0}^c (-1)^j \binom{[b]}j \times \binom{[a+c-j-1]}{a-1} \Longrightarrow \begin{cases} \binom{[a+c-b-1]}{c} & \mbox{ if } a \geq b\\ (-1)^c \binom{[b-a]}{c}  & \mbox{ if } a < b \end{cases}.$$
\end{problem}
\begin{proof}[Construction]
By the standard bijection between weak compositions\footnote{In weak compositions, $0$ is allowed as a summand.} of $n$ with $k$ summands and $\binom{[n+k-1]}{k-1}$, an element of the signed set on the left can be interpreted as a pair $(J,\pi)$, where $J$ is a subset of $[b]$, $\pi$ is a weak composition of $c-|J|$ with $a$ parts, and the sign of such a pair is $(-1)^{|J|}$.\\
 Assume first that $a \geq b$. The pairs with $J = \emptyset$ and $\pi_1 = \ldots = \pi_b = 0$ are clearly in bi\-jection with weak compositions of $c$ with $a - b$ parts, and therefore with $c$-element subsets of $[a+c-b-1]$. And we have a sign-reversing involution on the remaining pairs: if $i$ is the smallest element with $i \in J$ or $\pi_i > 0$, then remove $i$ from $J$ and add $1$ to $\pi_i$ if $i \in J$, and add $i$ to $J$ and subtract $1$ from $\pi_i$ if $i \notin J$.\\
 If $a < b$, the pairs with $J \cap [a] = \emptyset$ and $\pi = 0$ (and therefore necessarily $|J| = c$) have sign $(-1)^c$ and are clearly in bijection with subsets in $\binom{[b-a]}{c}$. On the remaining pairs, we have the same sign-reversing involution as before.
\end{proof}

If $a < b \leq a + c$, both sets $\binom{[a+c-b-1]}{c}$ and $\binom{[b-a]}{c}$ are empty. Note that the sijection proves the (obvious) equality
$$(1-x)^b \frac 1{(1-x)^a} = \begin{cases}\frac 1{(1-x)^{a-b}} & \mbox{if } a \geq b \\ (1-x)^{b-a} & \mbox{if } a < b \end{cases}.$$

\begin{problem} \label{prob:Bisijection}
 For $n,i \in \N$, construct a sijection
 $$\bigsqcup_{j=1}^n (-1)^{j+1} \binom{[2n-i-1]}{n-i-j+1} \times \B_{n,j} \Longrightarrow \B_{n,i}.$$
\end{problem}
\begin{proof}[Construction]
  By definition, we have a bijection $\B_{n,i} \to \binom{[n+i-2]}{n-1} \times \binom{[2n-i-1]}{n-1}$. So
  \begin{multline*}
    \bigsqcup_{j=1}^n (-1)^{j+1} \binom{[2n-i-1]}{n-i-j+1} \times \B_{n,j} \Rightarrow \bigsqcup_{j=1}^n (-1)^{j+1} \binom{[2n-i-1]}{n-i-j+1} \times \binom{[n+j-2]}{n-1} \times \binom{[2n-j-1]}{n-1} \\
    = \bigsqcup_{j=0}^{n-1} (-1)^{j} \binom{[2n-i-1]}{n-i-j} \times \binom{[n+j-1]}{n-1} \times \binom{[2n-j-2]}{n-1}.
  \end{multline*}
  By using the bijection from Problem \ref{prob:trinomial} for $a = n-i-j$, $b=n-1$, $c = j$, we get a sijection to
  $$\bigsqcup_{j=0}^{n-1} (-1)^{j} \binom{[2n-i-1]}{n-1} \times \binom{[n-i]}{j} \times \binom{[2n-j-2]}{n-1}.$$
  Now if use the sijection from Problem \ref{prob:binom} for $a = n$, $b = n-i$, $c = n-1$, we obtain a sijection to
  $$\binom{[2n-i-1]}{n-1} \times \binom{[n+i-2]}{n-1},$$
  which is equivalent to $\B_{n,i}$.
\end{proof}

The next goal is to prove that the sets $\ASM_{n,i}$ satisfy the same ``equalities'' as $\B_{n,i}$, i.e.~to construct a sijection
$$\bigsqcup_{j=1}^n (-1)^{j+1} \binom{[2n-i-1]}{n-i-j+1} \times \ASM_{n,j} \Longrightarrow \ASM_{n,i}.$$

An important part of the remainder of the construction will then be played by ``bijective linear algebra'' adapted to signed sets.

\section{Bijective linear algebra} \label{sect:LinAlg}

Denote by $\u{\mathfrak S}_m$ the signed set of permutations (with the usual sign). %More generally, given two (equipotent) finite subsets $A,B$ of $\Z$, denote by $\u{\mathfrak S}_{A,B}$ the signed set of all bijective maps $A \to B$, where the sign is $-1$ to the number of inversions.
Given signed sets $\u P_{i,j}$, $1 \leq i,j \leq m$, define the \emph{determinant} of $\u{\p P} = [\u P_{ij}]_{i,j = 1}^m$ as the signed set
$$\det(\u{\p P}) = \bigsqcup_{\pi \in \u{\mathfrak S}_m} \u P_{1,\pi(1)} \times \cdots \times \u P_{m,\pi(m)}.$$

Let us prove some classical properties.

\begin{problem} \label{prob:detproduct}
  Given $\u{\p P} = [\u P_{ij}]_{i,j = 1}^m$, $\u{\p Q} = [\u Q_{ij}]_{i,j = 1}^m$, construct a sijection
  $$\det(\u{\p R}) \Longrightarrow \det(\u{\p P}) \times \det(\u{\p Q}),$$
  where $\u{\p R} = [\u R_{ij}]_{i,j = 1}^m$, $\u R_{i,j} = \bigsqcup_{p=1}^m \u P_{i,p} \times \u Q_{p,j}$.
\end{problem}
\begin{proof}[Construction]
 On the left, we get (after using commutativity) signed sets
 $$\sign \pi \cdot \u P_{1, l_1} \times \cdots \times \u P_{m, l_m} \times \u Q_{l_1, \pi(1)} \times \cdots \times \u Q_{l_m, \pi(m)}$$
 for a permutation $\pi \in \u{\mathfrak S}_m$ and $l_1,\ldots,l_m \in [m]$. If there exist $i,j$, $i \neq j$, $l_i = l_j$ (and $i,j$ are the smallest such indices), then we have
 \begin{multline*}
   \sign \pi \cdot \u P_{1, l_1} \times \cdots \times \u P_{m, l_m} \times \u Q_{l_1, \pi(1)} \times \cdots \times \u Q_{l_m, \pi(m)} \\
   = - \sign \sigma \cdot \u P_{1, l_1} \times \cdots \times \u P_{m, l_m} \times \u Q_{l_1, \sigma(1)} \times \cdots \times \u Q_{l_m, \sigma(m)},
 \end{multline*}
 where $\sigma = \pi \cdot (i, j)$. If $(l_1,\ldots,l_m)$ form a permutation of $[m]$, then these are precisely the signed sets that appear in the Cartesian product $\det(\u{\p P}) \times \det(\u{\p Q})$, with the correct sign.
\end{proof}

\begin{problem} \label{prob:cramer}
 Given $\u{\p P} = [\u P_{p,q}]_{p,q = 1}^m$, signed sets $\u X_i, \u Y_i$  and sijections
 $$\bigsqcup_{q=1}^m \u P_{i,q} \times \u X_q \Longrightarrow \u Y_i$$
 for all $i \in [m]$, construct sijections
 $$ \det(\u{\p P}) \times \u X_j \Longrightarrow \det(\u{\p P}^j),$$
 where $\u{\p P}^j = [\u P^j_{p,q}]_{p,q = 1}^m$, $\u P^j_{p,q} = \u P_{p,q}$ if $q \neq j$, $\u P^j_{p,j} = \u Y_p$, for all $j \in [m]$.
\end{problem}
\begin{proof}[Construction]
 For a given permutation $\pi$ and fixed $i \in [m]$, we have sijections
 \begin{multline} \label{equ:det}
 \bigsqcup_{q = 1}^m \u P_{1, \pi(1)} \times \cdots \times  \u P_{i-1,\pi(i-1)} \times \u P_{i, q} \times \u P_{i+1,\pi(i+1)} \times \cdots \times \u P_{m,\pi(m)} \times \u X_q \\
 \Longrightarrow
  \u P_{1, \pi(1)} \times \cdots \times  \u P_{i-1,\pi(i-1)} \times \u Y_i \times \u P_{i+1,\pi(i+1)} \times \cdots \times \u P_{m,\pi(m)}
 \end{multline}
 by assumption. In the disjoint union over all permutations $\pi$ with $\pi(i)=j$ and over all $i \in [m]$ of the domains of these sijections, we obviously get the term $\det(\u{\p P}) \times \u X_j$ for $q=j$. If $q \neq j$, then
 \begin{multline*}
   \sign \pi \cdot \u P_{1,\pi(1)} \times \cdots \times \u P_{i-1,\pi(i-1)} \times \u P_{i,q} \times \u P_{i+1,\pi(i+1)} \times \cdots \times \u P_{m,\pi(m)} \\
   = - \left( \sign \sigma  \cdot \u P_{1,\sigma(1)} \times \cdots \times \u P_{p-1,\sigma(p-1)} \times \u P_{p,q} \times \u P_{p+1,\sigma(p+1)} \times \cdots \times \u P_{m,\sigma(m)}\right),
 \end{multline*}
  where $p$ is chosen so that $\pi(p) = q$, and $\sigma = \pi \cdot (i, p)$. In other words, there is a sijection from the coefficient at $\u X_q$ to $\u \emptyset$. On the other hand, the disjoint union over all permutations $\pi$ with $\pi(i)=j$ and over all $i \in [m]$ of the codomains of the sijections in \eqref{equ:det} is
just  $\det(\u{\p P}^j)$.
\end{proof}

The following is implied immediately, as $\det(\u{\p P}^j)$ is clearly $\u \emptyset$ if $\u Y_1 = \ldots = \u Y_m = \u \emptyset$.

\begin{problem} \label{prob:solveeq}
 Given $\u{\p P} = [\u P_{ij}]_{i,j = 1}^m$, signed sets $\u X_j$ for $j \in [m]$, and sijections
 $$\bigsqcup_{j=1}^m \u P_{i,j} \times \u X_j \Longrightarrow \u \emptyset,$$
 construct a sijection
 $$ \det(\u{\p P}) \times \u X_j \Longrightarrow \u \emptyset.$$
\end{problem}

\section{Prerequisites from Part I} \label{sect:pre}

In this section, we summarize definitions and results from Part I that are necessary for the constructions that will follow.

The following objects extend the classical Gelfand-Tsetlin patterns, see  \cite[p.\ 313]{Sta99} or \cite[(3)]{gelfand} for the original reference.

\begin{definition}
 For $k \in \Z$, define ${\GT(k) = (\{\cdot\},\emptyset)}$,
and for $\q k = (k_1,\ldots,k_n) \in \Z^n$, define {recursively}
 $$\GT(\q k) = \GT(k_1,\ldots,k_n) = \bigsqcup_{\q l \in \si{k_1}{k_2} \times \cdots \times \si{k_{n-1}}{k_n}} \GT(l_1,\ldots,l_{n-1}).$$
\end{definition}

In \cite{PartI}, we have constructed the following sijections.

\begin{problem}(\cite[Problem 5]{PartI}) \label{prob:pi}
  Given $\q k =(k_1,\ldots,k_n) \in \Z^n$ and $i$, $1 \leq i \leq n-1$, construct a sijection
 $$\pi = \pi_{\q k,i} \colon \GT(k_1,\ldots,k_n) \Rightarrow -\GT(k_1,\ldots,k_{i-1},k_{i+1}+1,k_i-1,k_{i+2},\ldots,k_n).$$
 \end{problem}

In the case that $k_1 \le k_2 \le \ldots \le k_n$, Gelfand-Tsetlin patterns can be thought of triangular arrays of integers of the following form
\begin{equation*}
\begin{array}{ccccccccccccccccc}
  &   &   &   &   &   &   &   & a_{1,1} &   &   &   &   &   &   &   & \\
  &   &   &   &   &   &   & a_{2,1} &   & a_{2,2} &   &   &   &   &   &   & \\
  &   &   &   &   &   & \dots &   & \dots &   & \dots &   &   &   &   &   & \\
  &   &   &   &   & a_{n-2,1} &   & \dots &   & \dots &   & a_{n-2,n-2} &   &   &   &   & \\
  &   &   &   & a_{n-1,1} &   & a_{n-1,2} &  &   \dots &   & \dots   &  & a_{n-1,n-1}  &   &   &   & \\
  &   &   & a_{n,1} &   & a_{n,2} &   & a_{n,3} &   & \dots &   & \dots &   & a_{n,n} &   &   &
\end{array},
\end{equation*}
such that $a_{n,i}=k_i$ for $i \in [n]$, and the entries increase weakly along $\nearrow$-diagonals and along $\searrow$-diagonals. In this classical case, we speak of a monotone triangle if the rows are strictly increasing. Monotone triangles with bottom row $(1,2,\ldots,n)$ are in easy bijective correspondence with ASMs of size $n \times n$, see, for instance, \cite{PartI}.

The notion of monotone triangles can also be extended to the case when
$\q k = (k_1,\ldots,k_n)$ is not weakly increasing as is explained next.

\begin{definition}
\begin{enumerate}
\item
Suppose that $\q k = (k_1,\ldots,k_n)$ and $\q l = (l_1,\ldots,l_{n-1})$ are two sequences of integers. We say that $\q l$ \emph{interlaces} $\q k$, $\q l \prec \q k$, if the following holds:
\begin{enumerate}
  \item for every $i$, $1 \leq i \leq n-1$, $l_i$ is in the closed interval between $k_i$ and $k_{i+1}$;
  \item if $k_{i-1} \leq k_i \leq k_{i+1}$ for some $i$, $2 \leq i \leq n-1$, then $l_{i-1}$ and $l_i$ cannot both be $k_i$;
  \item if $k_i > l_i = k_{i+1}$, then $i \leq n-2$ and $l_{i+1} = l_i = k_{i+1}$;
  \item if $k_i = l_i > k_{i+1}$, then $i \geq 2$ and $l_{i-1} = l_i = k_i$.
\end{enumerate}
\item A \emph{monotone triangle of size $n$} is a map $T \colon \{(i,j) \colon 1 \leq j \leq i \leq n \} \to \Z$ so that line $i-1$ (i.e.~the sequence $T_{i-1,1},\ldots,T_{i-1,i-1}$) interlaces line $i$ (i.e.~the sequence $T_{i,1},\ldots,T_{i,i}$).
\item The \emph{sign} of a monotone triangle $T$ is $(-1)^r$, where $r$ is the sum of:
\begin{itemize}
  \item the number of strict descents in the rows of $T$, i.e.~the number of pairs $(i,j)$ so that $1 \leq j < i \leq n$ and $T_{i,j} > T_{i,j+1}$, and
  \item the number of $(i,j)$ so that $1 \leq j \leq i - 2$, $i \leq n$ and $T_{i,j} > T_{i-1,j} = T_{i,j+1} = T_{i-1,j+1} > T_{i,j+2}$.
\end{itemize}
\end{enumerate}
\end{definition}

{The following is a monotone triangles of size $5$ with sign $-1$.}
$$\begin{array}{ccccccccc}
&&&& 4 &&&& \\
&&& 3 && 5 &&& \\
&& 3 && 4 && 5 && \\
& 3 & & 3 & & 4 & & 5 & \\
5 & & 3 & & 1 & & 4 & & 6
\end{array}$$

The purpose of \cite{PartI} was to construct a sijection between monotone triangles and yet another type of objects, namely shifted Gelfand-Tsetlin patterns, to be defined next. These new objects involve
 \emph{arrow patterns}, which are triangular arrays $T=(t_{p,q})_{1 \le p < q \le n}$ arranged as
$$T = \begin{smallmatrix} & & & & t_{1,n} & & & & \\ & & & t_{1,n-1} & & t_{2,n} & & & \\ & & t_{1,n-2} & & t_{2,n-1}& & t_{3,n} & & \\ & \vstretch{0.35}{\udots} & \vstretch{0.35}{\vdots} & \vstretch{0.35}{\ddots} & \vstretch{0.35}{\vdots} & \vstretch{0.35}{\udots} & \vstretch{0.35}{\vdots} & \vstretch{0.35}{\ddots} & \\ t_{1,2} & & t_{2,3}  & & \ldots  & & \ldots  & & t_{n-1,n}  \end{smallmatrix},$$
with $t_{p,q} \in \{\swarrow, \searrow, \seswarrow\}$. The sign of an arrow pattern is $1$ if the number of $\seswarrow$'s is even and $-1$ otherwise, and  the signed set of \emph{arrow patterns of order $n$} is denoted by $\AP_n$.

The role of an arrow pattern of order $n$ is that it induces a deformation of $(k_1,\ldots,k_n)$, which can be thought of as follows. Add $k_1,\ldots,k_n$ as bottom row of $T$ (i.e., $t_{i,i}=k_i$), and for each $\swarrow$ or$ \seswarrow$ which is in the same $\swarrow$-diagonal as $k_i$ add $1$ to $k_i$, while for each $\searrow$ or $\seswarrow$ which is in the same $\searrow$-diagonal as $k_i$ subtract $1$ from $k_i$.
More formally, letting $\delta_{\swarrow}(\swarrow) = \delta_{\swarrow}(\seswarrow) = \delta_{\searrow}(\searrow) = \delta_{\searrow}(\seswarrow) = 1$ and $\delta_{\swarrow}(\searrow) = \delta_{\searrow}(\swarrow) = 0$, we set
$$c_i(T) = \sum_{j=i+1}^{n} \delta_{\swarrow}(t_{i,j}) - \sum_{j=1}^{i-1} \delta_{\searrow}(t_{j,i}) \mbox{ and } d(\q k,T) = (k_1+c_1(T),k_2+c_2(T),\ldots,k_n+c_n(T))$$
for $\q k = (k_1,\ldots,k_n)$ and $T \in \AP_n$.

\begin{definition} For $\q k =(k_1,\ldots,k_n)$, define \emph{shifted Gelfand-Tsetlin patterns}, or SGT patterns for short, as the following disjoint union of GT patterns over arrow patterns of order $n$:
$$
\SGT(\q k) = \bigsqcup_{T \in \AP_n} \GT(d(\q k,T))
$$
\end{definition}

The main result of \cite{PartI} is a construction that solves the following problem.

\begin{problem}(\cite[Problem 10]{PartI})\label{gamma}
 Given $\q k = (k_1,\ldots,k_n) \in \Z^n$ and $x \in \Z$, construct a sijection
 $$\Gamma = \Gamma_{\q k,x} \colon \MT(\q k) \Rightarrow \SGT(\q k).$$
\end{problem}

A computer code for this construction is also available at the webpage mentioned in the introduction,
and the concrete example $\q k = (1,2,3)$ and $x=0$ is provided in \cite{PartI}.

\section{Rotation of monotone triangles}

The purpose of this section is solve the following problem.
\begin{problem} \label{prob:rot}
  Given $\q k = (k_1,\ldots,k_n)$, construct a sijection
  $$\MT(\q k) \Longrightarrow (-1)^{n-1} \MT(\rot(\q k)),$$
  where $\rot(\q k) = (k_2,\ldots,k_n,k_1-n)$.
\end{problem}

For a set $U \subseteq \N$, define $\epsilon_j(U) = 1$ if $j \in U$ and $0$ otherwise. For $\q k =(k_1,\ldots,k_n)$ and $j \geq 1$, define
$$\u E_j(\q k) = \bigsqcup_{V \in \binom{[n]}j} \bigsqcup_{U \subseteq V} (-1)^{|U|} \GT(k_1+\epsilon_1(U),k_2+\epsilon_2(U),\ldots,k_n+\epsilon_n(U))$$
and
$$\u F_j(\q k) = \bigsqcup_{V \in \binom{[n]}j} \bigsqcup_{U \subseteq V} (-1)^{|U|} \GT(k_1-\epsilon_1(U),k_2-\epsilon_2(U),\ldots,k_n-\epsilon_n(U)).$$

\begin{problem} \label{prob:prepelementary}
Given $\q k = (k_1,\ldots,k_n)$ and $j \geq 0$, construct a normal sijection
\begin{multline*}
\bigsqcup_{V \in \binom{[n]}j} \bigsqcup_{U \subseteq V} (-1)^{|U|}
\si{k_1+\epsilon_1(U)}{k_2+\epsilon_2(U)} \times \cdots \times \si{k_{n-1}+\epsilon_{n-1}(U)}{k_{n}+\epsilon_n(U)}  \\
\Longrightarrow \bigsqcup_{V \in \binom{[n-1]}j} \bigsqcup_{U \subseteq V} \bigsqcup_{\substack{I \subseteq [n-2] \\ I \cap (I + 1) = \emptyset}} (-1)^{|U|} \prod_{i=1}^{n-1} \begin{cases}
\{k_{i+1}+1\}  & i \in I \\
\{ k_i \} & i-1 \in I \\
\si{k_i+\epsilon_i(U)}{k_{i+1}+\epsilon_i(U)} & i-1,i \notin I
\end{cases}.
\end{multline*}
\end{problem}

\begin{proof}[Construction]
Construction is by induction with respect to $n$. Consider $n=2$. Both sides are empty if $j>2$, and both sides are also obviously equal if $j=0$. The cases $j=1,2$ follow easily from the construction in Problem~\ref{prop:alpha}.

Suppose $n \ge 3$. Then there exist normal sijections
\begin{multline*}
\bigsqcup_{V \in \binom{[n]}j} \bigsqcup_{U \subseteq V} (-1)^{|U|} \prod_{i=1}^{n-1}
\si{k_i+\epsilon_i(U)}{k_{i+1}+\epsilon_{i+1}(U)}  \\
\Longrightarrow \bigsqcup_{V \in \binom{[n-1]}j} \bigsqcup_{U \subseteq V} (-1)^{|U|} \left( \prod_{i=1}^{n-2}
\si{k_i+\epsilon_i(U)}{k_{i+1}+\epsilon_{i+1}(U)} \right) \times \si{k_{n-1}+\epsilon_{n-1}(U)}{k_{n}} \\ \sqcup
 \bigsqcup_{V \in \binom{[n-1]}{j-1}} \bigsqcup_{U \subseteq V} (-1)^{|U|+1}
 \left( \prod_{i=1}^{n-2}
\si{k_i+\epsilon_i(U)}{k_{i+1}+\epsilon_{i+1}(U)} \right) \times \si{k_{n-1}+\epsilon_{n-1}(U)}{k_{n}+1} \\
\Longrightarrow \bigsqcup_{V \in \binom{[n-1]}j} \bigsqcup_{U \subseteq V} (-1)^{|U|}  \left( \prod_{i=1}^{n-2}
\si{k_i+\epsilon_i(U)}{k_{i+1}+\epsilon_{i+1}(U)} \right) \times  \si{k_{n-1}}{k_{n}} \\ \sqcup \bigsqcup_{V \in \binom{[n-1]}j} \bigsqcup_{U \subseteq V} (-1)^{|U|+1} \left( \prod_{i=1}^{n-2}
\si{k_i+\epsilon_i(U)}{k_{i+1}+\epsilon_{i+1}(U)} \right) \times \si{k_{n-1}}{k_{n-1}+\epsilon_{n-1}(U)-1} \\
\sqcup
 \bigsqcup_{V \in \binom{[n-1]}{j-1}} \bigsqcup_{U \subseteq V} (-1)^{|U|+1}
 \left( \prod_{i=1}^{n-2}
\si{k_i+\epsilon_i(U)}{k_{i+1}+\epsilon_{i+1}(U)} \right) \times \si{k_{n-1}+1}{k_{n}+1} \\
\sqcup
 \bigsqcup_{V \in \binom{[n-1]}{j-1}} \bigsqcup_{U \subseteq V} (-1)^{|U|}
 \left( \prod_{i=1}^{n-2}
\si{k_i+\epsilon_i(U)}{k_{i+1}+\epsilon_{i+1}(U)} \right) \times \si{k_{n-1}+1}{k_{n-1}+\epsilon_{n-1}(U)-1},
\end{multline*}
using $\alpha$ from Problem~\ref{prop:alpha} twice as well as $\si{a}{b} = - \si{b+1}{a-1}$. By induction, the first and the third term can be combined to obtain the terms in the codomain with $n-2 \notin I$.

The second term is equal to $\u \emptyset$ unless $\epsilon_{n-1}(U)=1$, while the fourth term is equal to $\u \emptyset$ unless $\epsilon_{n-1}(U)=0$. Combining these terms we get
\begin{multline*}
\bigsqcup_{V \in \binom{[n-2]}{j-1}} \bigsqcup_{U \subseteq V} (-1)^{|U|} \left( \prod_{i=1}^{n-3}
\si{k_i+\epsilon_i(U)}{k_{i+1}+\epsilon_{i+1}(U)} \right) \times
\si{k_{n-2}+\epsilon_{n-2}(U)}{k_{n-1}+1} \times \u {\{ k_{n-1} \}} \\
\sqcup \bigsqcup_{V \in \binom{[n-2]}{j-1}} \bigsqcup_{U \subseteq V} (-1)^{|U|+1} \left( \prod_{i=1}^{n-3}
\si{k_i+\epsilon_i(U)}{k_{i+1}+\epsilon_{i+1}(U)} \right) \times \si{k_{n-2}+\epsilon_{n-2}(U)}{k_{n-1}}  \times \u {\{ k_{n-1}\} }  \\
\Longrightarrow \bigsqcup_{V \in \binom{[n-2]}{j-1}} \bigsqcup_{U \subseteq V} (-1)^{|U|} \left( \prod_{i=1}^{n-3}
\si{k_i+\epsilon_i(U)}{k_{i+1}+\epsilon_{i+1}(U)} \right) \times \u {\{ k_{n-1}+1 \}} \times \u {\{ k_{n-1} \}},
\end{multline*}
using $\alpha$ to obtain the sijection. By induction, the codomain of this sijection is the term of the codomain in the problem with $n-2 \in I$.

The combination of the sijections is normal because $\alpha$ is normal.
\end{proof}

\begin{problem} \label{prob:elementary}
  Given $\q k = (k_1,\ldots,k_n)$ and $j \geq 1$, construct sijections
  $$\u E_j(\q k) \Longrightarrow \u \emptyset \quad \mbox{and} \quad \u F_j(\q k) \Longrightarrow \u \emptyset.$$
\end{problem}
\begin{proof}[Construction]
 It is clear that $\u E_j(\q k) = \u \emptyset$ if $j > n$, so it is enough to construct a sijection
 $$\u E_j(\q k) \Longrightarrow \bigsqcup_{\q l \in \si{k_1}{k_2} \times \cdots \times \si{k_{n-1}}{k_n}} \u E_j(\q l).$$ We obtain such a sijection using the normal sijections from
 Problem~\ref{prob:prepelementary}: Observe that $\u E_j(\q k)$ is
 the disjoint union of the family $\GT(l_1,\ldots,l_{n-1})$ over the domain of the sijection in Problem~\ref{prob:prepelementary}, while $\bigsqcup_{\q l \in \si{k_1}{k_2} \times \cdots \times \si{k_{n-1}}{k_n}} \u E_j(\q l)$ is the disjoint union of the same family over the part of the codomain where $I=\emptyset$. For the part of the codomain where $I \not= \emptyset$,
 $\GT(l_1,\ldots,l_{n-1}) = \u \emptyset$ as $\GT(l_1,\ldots,l_{n-1})  = \u \emptyset$ whenever $l_{i}=l_{i+1}+1$ for some $i$ since $\si{x+1}{x} = \u \emptyset$. The sijection is then the disjoint union of the identity sijections as discussed at the end of
 Section~\ref{sect:signedsets}.

The construction of the sijection $\u F_j(\q k) \Longrightarrow \u \emptyset$ is analogous.
\end{proof}

Let us also define

$$\u E_j'(\q k) = \bigsqcup_{V \in \binom{[2,n]}j} \bigsqcup_{U \subseteq V} (-1)^{|U|} \GT(k_1,k_2+\epsilon_2(U),\ldots,k_n+\epsilon_n(U))$$
and
$$\u F_j'(\q k) = \bigsqcup_{V \in \binom{[2,n]}j} \bigsqcup_{U \subseteq V} (-1)^{|U|} \GT(k_1,k_2-\epsilon_2(U),\ldots,k_n-\epsilon_n(U)).$$

The previous problem can be rephrased as follows.

\begin{problem} \label{prob:elementaryprime}
  Given $\q k = (k_1,\ldots,k_n)$ and $j \geq 0$, construct sijections
  $$\u E_j'(\q k) \Longrightarrow \bigsqcup_{i = 0}^j (-1)^{j-i} \binom{[j]}i \times \GT(k_1+i,k_2,\ldots,k_n)$$
  and
  $$\u F_j'(\q k) \Longrightarrow \bigsqcup_{i = 0}^j (-1)^{j-i}  \binom{[j]}i \times \GT(k_1-i,k_2,\ldots,k_n).$$
\end{problem}
\begin{proof}[Construction]
 The sijections are trivial for $j = 0$. For $j > 0$, we can split pairs $(U,V)$, $U \subseteq V \in \binom{[n]}j$, that appear in $\u E_j(\q k)$ into three groups: $1 \notin V$, $1 \in U$, and $1 \in V, 1 \notin U$. This gives an equivalence $\u E_j(\q k) \approx \u E_j'(\q k) \sqcup (-\u E_{j-1}'(k_1+1,k_2,\ldots,k_n)) \sqcup \u E_{j-1}'(\q k)$ and, together with Problem \ref{prob:elementary}, an equivalence $\u E_j'(\q k) \approx \u E_{j-1}'(k_1+1,k_2,\ldots,k_n) \sqcup (- \u E_{j-1}'(\q k))$. The first sijection is now constructed inductively, and the second construction is analogous.
\end{proof}

The following problem involves arrow patterns (see Section~\ref{sect:pre}) as well as arrow rows, which are defined as follows: An \emph{arrow row of order $n$} is a row of length $n$ with elements $\nwarrow, \nearrow, \nenwarrow$, where the positive elements are precisely those with an even number of $\nenwarrow$'s. The signed set of arrow rows of length $n$ is denoted by $\AR_n$. Given $T' \in \AP_{n-1}$ and $\mu \in \AR_n$, write $\mu T'$ (resp., $T'\mu$) for the arrow pattern of size $n$ that we obtain if we add $\mu$ (with all arrows reflected along the horizontal axis) to $T'$ as the leftmost $\nearrow$-diagonal (resp., rightmost $\searrow$-diagonal) in such a way that the first element of the arrow row is the bottom (resp., top) element of the added diagonal.

\begin{problem} \label{prob:lefttoright}
  Given $\q k= (k_1,\ldots,k_n)$ and $T' \in \AP_{n-1}$, construct a sijection
  \begin{multline*}
   \bigsqcup_{\mu \in \AR_{n-1}} \GT(k_1+c_1(\mu T'),k_2+c_2(\mu T'),\ldots,k_n+c_{n}(\mu T')) \\
   \Longrightarrow \bigsqcup_{\mu \in \AR_{n-1}} \GT(k_1+c_n(T'\mu),k_2+c_1(T'\mu),\ldots,k_n+c_{n-1}(T'\mu)).
  \end{multline*}
\end{problem}
\begin{proof}[Construction]
 Write $l_i = k_i + c_{i-1}(T')$ for $i = 2,\ldots,n$, and define
 $$\u F = \bigsqcup_{j = 0}^{n-1} (-1)^{n-1-j} \bigsqcup_{p = 0}^{n-1-j} (-1)^p \binom{[n-j-1]}p \times \u F_{n-j-1}'(k_1+p,l_2,\ldots,l_n).$$
 By the definition of $\u F_j'$, this signed set is equivalent to
 $$\bigsqcup_{p, \epsilon} (-1)^p \left( \bigsqcup_{j = 0}^{n-1-p}  \bigsqcup_{V} (-1)^{n - 1 - j + q} \binom{[n-j-1]}p \right) \times \GT(k_1+p, l_2 - \epsilon_2,\ldots,l_n - \epsilon_n),$$
 where the outer disjoint union is over $p \in \si 0{n-1}$ and $\epsilon = (\epsilon_2,\ldots,\epsilon_{n}) \in \si 0 1 ^n$, the innermost disjoint union is over $V \in \binom{[2,n]}{n-j-1}$ that contain all $i$ with $\epsilon_i = 1$, and we set $q = \epsilon_2 + \ldots + \epsilon_n$. Clearly, this is equivalent to
 $$\bigsqcup_{p, \epsilon} (-1)^p \left( \bigsqcup_{j = 0}^{n-1-p}  (-1)^{n - 1 - j + q} \binom{[n-1-q]}{n-1-q-j} \times \binom{[n-j-1]}p \right) \times \GT(k_1+p, l_2 - \epsilon_2,\ldots,l_n - \epsilon_n),$$
 and Problem \ref{prob:binom} for $a = p + 1$, $b = n - 1 - q$ and $c = n - p - 1$ gives us the equivalence
 $$\u F \approx \bigsqcup_{p, \epsilon} (-1)^{n + p + q - 1} \binom{[q]}{n-p-1} \times \GT(k_1+p, l_2 - \epsilon_2,\ldots,l_n - \epsilon_n).$$
 On the other hand, note that $k_1 + c_1(\mu T') = k_1 + \sum_{i=1}^{n-1} \delta_{\nwarrow}(\mu_i)$ and $k_i + c_i(\mu T') = l_i - \delta_{\nearrow}(\mu_{i-1})$ for $i = 2,\ldots,n$. Let us see how many times $\GT(k_1 + p, l_2 - \epsilon_2,\ldots,l_n - \epsilon_n)$ appears in the disjoint union $\bigsqcup_{\mu \in \AR_{n-1}} \GT(k_1+c_1(\mu T'),k_2+c_2(\mu T'),\ldots,k_n+c_{n}(\mu T'))$. Whenever $\epsilon_i = 1$, we must have $\mu_{i-1} \in \{\nearrow, \nenwarrow\}$, and whenever $\epsilon_i = 0$, we must have $\mu_{i-1} = \nwarrow$. That means that among the $q$ $i$'s with $\epsilon_i = 1$, we have to select $p - (n-1-q)$ $\nenwarrow$'s. That gives an equivalence
 $$\u F \approx \bigsqcup_{\mu \in \AR_{n-1}} \GT(k_1+c_1(\mu T'),k_2+c_2(\mu T'),\ldots,k_n+c_{n}(\mu T')).$$
 In an analogous way, we prove that
 $$\u E \approx \bigsqcup_{\mu \in \AR_{n-1}} \GT(k_1+c_n(T'\mu),k_2+c_1(T'\mu),\ldots,k_n+c_{n-1}(T'\mu)),$$
 where
 $$\u E = \bigsqcup_{j = 0}^{n-1} (-1)^{n-1-j} \bigsqcup_{p = 0}^{n-1-j} (-1)^p \binom{[n-j-1]}p \times \u E_{n-j-1}'(k_1-p,l_2,\ldots,l_n).$$
 The sijection $\u E \Rightarrow \u F$ easily follows from Problem \ref{prob:elementaryprime}.
\end{proof}

Now we are in the position to give a sijective proof of the rotational invariance up to sign.

\begin{proof}[Construction for Problem~\ref{prob:rot}]
We use the sijection $\Gamma$ from Problem~\ref{gamma} (which depends on a parameter $x$).
 \begin{multline*}
  \MT(\q k) \Longrightarrow \SGT(\q k) = \bigsqcup_{T \in \AP_n} \GT(k_1 + c_1(T), \ldots, k_n + c_n(T)) \\
 = \bigsqcup_{T' \in \AP_{n-1}} \bigsqcup_{\mu \in \AR_{n-1}} \GT(k_1+c_1(\mu T'),k_2+c_2(\mu T'),\ldots,k_n+c_{n}(\mu T')).
\end{multline*}
Apply the sijection from Problem \ref{prob:lefttoright}, and then rotation of GT patterns (using the sijections from Problem~\ref{prob:pi}). We get
\begin{multline*}
  \bigsqcup_{T' \in \AP_{n-1}} \bigsqcup_{\mu \in \AR_{n-1}} \GT(k_1+c_n(T'\mu),k_2+c_1(T'\mu),\ldots,k_n+c_{n-1}(T'\mu)) \\
  \Longrightarrow (-1)^{n-1} \bigsqcup_{T' \in \AP_{n-1}} \bigsqcup_{\mu \in \AR_{n-1}} \GT(k_2+c_1(T'\mu)+1,\ldots,k_n+c_{n-1}(T'\mu)+1,k_1+c_n(T'\mu)-n+1) \\
  = (-1)^{n-1} \bigsqcup_{T \in \AP_n} \GT(k_2 + c_1(T)+1, \ldots, k_n  +c_{n-1}(T)+1,k_1 + c_n(T)-n+1).
\end{multline*}
We subtract $1$ from all GT patterns and apply $\Gamma^{-1}$, and obtain $(-1)^{n-1} \MT(\rot(\q k))$.
\end{proof}

\section{Equalities determining $\ASM_{n,i}$}

Suppose that we are given a weakly increasing sequence $\q k = (k_1,\ldots,k_n)$ and $i \in \N$. We define
$$\MT_i(\q k) = \{ T \in \MT(\q k) \colon T_{n-i+1,1} = \ldots = T_{n,1} = k_1, T_{n-i,1} \neq k_1\}$$
as the signed subset of monotone triangles with $k_1$ in the first position in exactly the last $i$ rows. Similarly, we define
$$\MT^i(\q k) = \{ T \in \MT(\q k) \colon T_{n-i+1,n-i+1} = \ldots = T_{n,n} = k_n, T_{n-i,n-i} \neq k_n\}$$
as the signed subset of monotone triangles with $k_n$ in the last position in exactly the last $i$ rows.

\begin{problem} \label{prob:mti}
 Given a weakly increasing $\q k = (k_1,\ldots,k_n)$ and $i \geq 1$, construct sijections
 $$\MT_i(\q k) \Longrightarrow \bigsqcup_{j = 0}^{i-1} (-1)^j \binom{[i-1]}j \times \MT(k_1+j+1,k_2,\ldots,k_n)$$
 and
 $$\MT^i(\q k) \Longrightarrow \bigsqcup_{j = 0}^{i-1} (-1)^j \binom{[i-1]}j \times \MT(k_1,k_2,\ldots,k_n-j-1).$$
\end{problem}
\begin{proof}[Construction]
 The signed sets $\MT_1(\q k)$ and $\MT(k_1+1,k_2,\ldots,k_n)$ are obviously equivalent: just increase/decrease the first element of the bottom row of a monotone triangle to get a sijection. For $i > 1$, we have  $$\MT_i(\q k) = \bigsqcup_{\q l} \MT_{i-1}(k_1,l_2,\ldots,l_{n-1}),$$
 where the disjoint union is over all $\q l = (l_2,\ldots,l_{n-1})$ that interlace $(k_2,\ldots,k_n)$. By induction, we have a sijection to
 $$\bigsqcup_{\q l} \bigsqcup_{j = 0}^{i-2} (-1)^j \binom{[i-2]}j \times \MT(k_1+j+1,l_2,\ldots,l_n),$$
 and it remains to find a sijection from this signed set to
 $$\bigsqcup_{j = 0}^{i-1} (-1)^j \binom{[i-1]}j \times \MT(k_1+j+1,k_2,\ldots,k_n).$$
 An element of this signed set is $((A,T),j)$ for a subset $A$ of $[i-1]$, $j = |A|$, and $T$ a monotone triangle with bottom row $(k_1+j+1,k_2,\ldots,k_n)$. Let the second-to-last row be $(k_1+p+1,l_2,\ldots,l_{n-1})$, where $(l_2,\ldots,l_{n-1})$ interlaces $(k_2,\ldots,k_n)$. If $k_1 + p + 1 < k_2$ or $k_1 + p + 1 = k_2 < l_2$, we must have $k_1 + j + 1 \leq k_1 + p + 1$, so we have two options. If $|A| = p$ and $i-1 \notin A$, the element $(((A,T),p),\q l)$ appears in the signed set on the left. The elements with $|A| < p$ or $i - 1 \in A$ cancel each other: simply add or remove $i-1$ from $A$. On the other hand, if $k_1 + p +1 > k_2$ or $k_1+p + 1 = k_2 = l_2$, we have $k_1 + j + 1 > k_1 + p  + 1$, and we again have two options. If $|A| = p + 1$ and $i-1 \in A$, the element $(((A \setminus \{i-1\},T),p),\q l)$ appears in the signed set on the left, and the elements with $|A| > p + 1$ or $i - 1 \notin A$ cancel each other.
\end{proof}

\begin{problem}  \label{prob:Aisijection}
  Given $n \in \N$ and $i \in [n]$, construct a sijection
  $$\bigsqcup_{j=1}^n (-1)^{j+1} \binom{[2n-i-1]}{n-i-j+1} \times \ASM_{n,j} \Longrightarrow \ASM_{n,i}. $$
\end{problem}
\begin{proof}[Construction]
 Using the obvious bijection $\ASM_{n,i} \to \ASM_{n,n+1-i}$ obtained by reflecting along the vertical axis and complementation, we get a sijection
 \begin{multline*}
 \bigsqcup_{j=1}^n (-1)^{j+1} \binom{[2n-i-1]}{n-i-j+1} \times \ASM_{n,j} \Longrightarrow \bigsqcup_{j=1}^n (-1)^{j+1} \binom{[2n-i-1]}{n-i-j+1} \times \ASM_{n,n+1-j} \\
 = \bigsqcup_{j=1}^\infty (-1)^{n-j} \binom{[2n-i-1]}{j-i} \times \ASM_{n,j}  \Longrightarrow \bigsqcup_j (-1)^{n-j} \binom{[2n-i-1]}{2n-j-1} \times \ASM_{n,j}\\
  = \bigsqcup_j (-1)^{n-j-1} \binom{[2n-i-1]}{j} \times \ASM_{n,2n-j-1}.
\end{multline*}
 By rotating the alternating sign matrix by $90^\circ$ counterclockwise and using the usual bijection between ASMs and monotone triangles with bottom row $1,\ldots,n$, we get a sijection $\ASM_{n,i} \to \MT_i(1,\ldots,n)$, and by reflecting along the vertical axis of symmetry and rotation by $90^\circ$ clockwise, we get a sijection $\ASM_{n,i} \to \MT^i(1,\ldots,n)$. The first sijection constructed in Problem \ref{prob:mti} gives us a sijection
 \begin{multline*}
 \bigsqcup_j (-1)^{n-j-1} \binom{[2n-i-1]}{j} \times \ASM_{n,2n-j-1} \\
 \Longrightarrow \bigsqcup_j (-1)^{n-j-1} \binom{[2n-i-1]}{j} \times \left( \bigsqcup_p (-1)^p\binom{[2n-j-2]}p \times \MT(2+p,2,\ldots,n)\right) \\
 \Longrightarrow  \bigsqcup_{p} (-1)^{n+p-1} \left( \bigsqcup_j (-1)^{j} \binom{[2n-i-1]}{j} \times \binom{[2n-j-2]}p \right)\times \MT(2+p,2,\ldots,n).
\end{multline*}
 The sijection from Problem \ref{prob:binom} for $a = p+1$, $b=2n-i-1$, $c=2n-p-2$ gives a sijection to
 $$\bigsqcup_{p} (-1)^{n+p-1} \binom{[i-1]}{2n-p-2} \times \MT(2+p,2,\ldots,n) = \bigsqcup_{p} (-1)^{n+p-1} \binom{[i-1]}{p} \times \MT(2n-p,2,\ldots,n).$$
 Now note that $\rot(2n-p,2,\ldots,n) = (2,\ldots,n,n-p)$ and use the sijection from Problem \ref{prob:rot} to get a sijection to
 $$\bigsqcup_{p} (-1)^p \binom{[i-1]}{p} \times \MT(2,\ldots,n,n-p) \Longrightarrow \bigsqcup_{p} (-1)^p \binom{[i-1]}{p} \times \MT(1,\ldots,n-1,n-1-p).$$
 The second sijection from Problem \ref{prob:mti} now gives a sijection to $\MT^i(1,\ldots,n) \approx \ASM_{n,i}$.
\end{proof}

\begin{problem} \label{prob:fromdet}
  Given $n$, construct a sijection
  $$\det(\u{\p P}) \Longrightarrow (-1)^{n-1} \DPP_{n-1},$$
  where $\u{\p P} = [\u P_{i,j}]_{i,j=2}^{n}$ with
$$\u P_{i,j} = (-1)^{j+1} \binom{[2n-i-1]}{n-i-j+1} \sqcup \begin{cases}
-\si 0 0 & i=j \\
\u \emptyset & \text{otherwise}
\end{cases}.
$$
%  for $2 \leq i,j \leq n$, $i \neq j$, and $\u P_{i,i} = (-1)^{i+1} \binom{[2n-i-1]}{n-2i+1} \sqcup -\si 0 0$ for $i=2,\ldots,n$.
\end{problem}
\begin{proof}[Construction]
% For $i = n$, we have $\binom{[2n-i-1]}{n-i-j+1} = \emptyset$ for $j = 2,\ldots,n$, so $\det(\u{\p P}) = - \det([\u P_{i,j}]_{i,j=2}^{n-1})$, and we set
%$\u {\p P'} = [\u P_{i,j}]_{i,j=2}^{n-1}$.
We define signed sets $\u S_{i,j} = (-1)^{i+j} \binom{[n]}{j-i}$ and $\u {\p S} = [\u S_{i,j}]_{i,j=1}^{n-1}$. It is clear that $\det(\u{\p S})$ is a signed set with one positive element and no negative elements, so, by Problem \ref{prob:detproduct}, we have a sijection
 $$\det(\u{\p P}) \Longrightarrow \det(\u{\p S})  \times \det(\u{\p P})  \Longrightarrow \det(\u{\p R}),$$
 where $\u{\p R} = [\u R_{ij}]_{i,j = 1}^{n-1}$, $\u R_{i,j} = \bigsqcup_{p=1}^{n-1} \u S_{i,p} \times \u P_{p+1,j+1}$. By using Problem \ref{prob:binom} for $a=n+j, b=n, c=n-i-j-1$,
\begin{multline*}
\u R_{i j} = \bigsqcup_{p=1}^{n-1} (-1)^{i+j+p} \binom{[n]}{p-i} \binom{[2n-p-2]}{n-p-j-1}
\longrightarrow  \bigsqcup_{p=1}^{n-1} (-1)^{i+j+p} \binom{[n]}{p-i} \binom{[2n-p-2]}{n+j-1}  \\
 = \bigsqcup_{k=0}^{n} (-1)^{k+j} \binom{[n]}{k} \binom{[2n-i-k-2]}{n+j-1}
\Longrightarrow  (-1)^j \binom{[n-i-2]}{n-i-j-1} \longrightarrow (-1)^j \binom{[n-i-2]}{j-1},
\end{multline*}
we get a sijection to $\det(\u{\p T})$, where $\u T_{i,j} = (-1)^j \binom{[n-i-2]}{j-1} \sqcup (-1)^{i+j+1} \binom{[n]}{j-i}$ for $1 \leq i,j \leq n -1$.\\
 We can also start with a matrix $\u{\p W} = [\u W_{i,j}]_{i,j=1}^{n-1}$, $\u W_{i,j} = \binom{[i+j]}{j-1}$ if $i \neq j$, and $\u W_{i,i} = \binom{[2i]}{i-1} \sqcup \si 0 0$. Its determinant can be, by the Lindstr\"om-Gessel-Viennot lemma \cite{Lin73,GesVie85}, interpreted as a collection of non-intersecting lattice paths between $(0,i+1)$ and $(i-1,0)$ for $i$ in a subset of $[n-2]$. By adding the step $(-1,i+1) \to (0,i+1)$, ignoring the steps at height $0$, and recording the heights of all steps, we indeed get descending plane partitions with parts $\leq n-1$. Again, we have sijections
 $$\det(\u{\p W}) \Longrightarrow  \det(\u{\p W}) \times  \det(\u{\p S})  \Longrightarrow \det(\u{\p U}),$$
 where $\u U_{i,j} = \bigsqcup_{p=1}^{n-1} \u W_{i,p} \times \u S_{p,j}$. Now we use Problem \ref{prob:binom} for $a = i+2$, $b = n$, $c = j-1$,
\begin{multline*}
\u U_{i,j} = \bigsqcup_{p=1}^{n-1}  (-1)^{p+j} \binom{[i+p]}{p-1} \binom{[n]}{j-p} \longrightarrow
\bigsqcup_{p=1}^{n-1}  (-1)^{p+j} \binom{[i+p]}{i+1} \binom{[n]}{j-p} \\
= \bigsqcup_{k=0}^{n} (-1)^k \binom{[i+j-k]}{i+1} \binom{[n]}{k} \Longrightarrow (-1)^{j-1} \binom{[n-i-2]}{j-1},
 \end{multline*}
which gives a sijection to $(-1)^{n-1}\det(\u{\p T})$. This finishes the construction.
\end{proof}

Now we have to put all the ingredients in place.

\begin{proof}[Construction for Problem \ref{prob:main}]
By Problems \ref{prob:solveeq} and \ref{prob:fromdet}, it is enough to construct sijections
$$\bigsqcup_{j = 2}^n \u P_{i,j} \times \left( \B_{n,1} \times \ASM_{n,j} \sqcup - \ASM_{n,1} \times \B_{n,j}\right) \Rightarrow \u \emptyset,$$
where $\u P_{i,j}$ was defined in Problem \ref{prob:fromdet}. But this is done by adding the signed set $\u P_{i,1} \times (\B_{n,1} \times \ASM_{n,1} \sqcup - \ASM_{n,1} \times \B_{n,1})$, which has size $0$, and using Problems \ref{prob:Bisijection} and \ref{prob:Aisijection} on the two terms of the resulting disjoint union.
\end{proof}

\begin{proof}[Construction for Problem~\ref{prob:asmtodpp}]
Using Problem~\ref{prob:Aisijection}, there are sijections
$$
\bigsqcup_{j=2}^n - \u P_{i,j}  \times \ASM_{n,j}  \Longrightarrow \binom{[2n-i-1]}{n-i} \times \ASM_{n,1}
$$
for all $i=1,2,\ldots,n$, where $\u P_{i,j}$ is defined as in Problem~\ref{prob:fromdet}. Using Problem~\ref{prob:cramer} and Problem~\ref{prob:fromdet}, it follows that there are sijections
$$
\DPP_{n-1} \times \ASM_{n,j} \Longrightarrow \det( \u P^j),
$$
where $\u P^j$ is obtained from $\u {\p P}$ by replacing the $j$-th column by
$\left( \binom{[2n-i-1]}{n-i} \times \ASM_{n,1} \right)^T_{2 \le i \le n}$. Equivalently, with
$\u Q_{i} =  \binom{[2n-i-1]}{n-i} \times \ASM_{n,1}$
and  $\u {\p Q}^j$ the $(n-1) \times (n-1)$ matrix of all empty sets, except for the $j$-th column, $2 \le j \le n$, which is $[\u Q_i]_{i=2}^{n}$, there is a sijection
$$
\DPP_{n-1} \times \ASM_{n,j} \Longrightarrow \det( \u {\p P} \sqcup \u {\p Q}^j) \sqcup -\det( \u {\p P}).
$$
Using $\u {\p S}$ from Problem~\ref{prob:fromdet}, we obtain
$$
\det( \u {\p P} \sqcup \u {\p Q}^j) \sqcup -\det( \u {\p P}) \Longrightarrow \det( \u {\p S}) \times \left( \det( \u {\p P} \sqcup \u {\p Q}^j) \sqcup -\det( \u {\p P}) \right).
$$
Using Problem \ref{prob:binom} for $a=n, b=n, c=n-i-1$,
 \begin{multline*} \bigsqcup_{p=1}^{n-1} \u S_{i,p} \times \u Q_{p+1} =
 \bigsqcup_{p=1}^{n-1} (-1)^{i+p}  \binom{[n]}{p-i} \times \binom{[2n-p-2]}{n-1} \times \ASM_{n,1} \\
= \bigsqcup_{k=0}^{n} (-1)^{k}  \binom{[n]}{k} \times \binom{[2n-i-k-2]}{n-1} \times \ASM_{n,1} \Longrightarrow
\begin{cases} \ASM_{n,1} & i=n-1 \\ \u \emptyset & i<n-1 \end{cases}.
\end{multline*}
%Letting $\u {\p Q}'^j=[\u Q_{p,q}]_{p,q,=1}^{n-1}$ with $\u Q_{p,q}
%be the $(n-1) \times (n-1)$ matrix of all empty sets, except for the entry in the $(j-1)$-st column, $2 \le j \le n$, and last row, which is $\ASM_{n,1}$.
Using $\u {\p T}$ from Problem~\ref{prob:binom}, we obtain a sijection to
$$
%\det( \u {\p T} \sqcup \u {\p Q}'^j) \sqcup -\det( \u {\p T}) \Longrightarrow
(-1)^{n+j} \det( \u {\p T}^j) \times \ASM_{n,1},
$$
where $\u {\p T}^j$ is obtain from $\u {\p T}$ by deleting the last row and $(j-1)$-st column.

On the other hand, let $\u {\p W}^j = [ \u W_{p,q} ]_{i,j=1}^{n-1}$ with
$$
\u W_{p,q} = \begin{cases} \binom{[p+q]}{q-1} & p < n-1 \\ \binom{[n+q-j]}{n-1} & p=n-1\end{cases} \sqcup
\begin{cases} [0,0] & p=q < n-1 \\ \u \emptyset & \text{otherwise} \end{cases}.
$$
By refining the considerations from Problem~\ref{prob:fromdet}, there is a sijection from $\det(\u{\p W}^j)$ to $\DPP_{n,j}$ using the Lindstr\"om-Gessel-Viennot lemma.
Again we have sijections
$$\det(\u{\p W}^j) \Longrightarrow  \det(\u{\p W}^j) \times  \det(\u{\p S})  \Longrightarrow \det(\u{\p U}^j),$$
with $\u{\p U}^j = [\u U_{p,q}]_{p,q=1}^{n-1}$ such that $\u{\p U}^j$ coincides on the top $n-2$ rows to $\u{\p T}$ up to sijections. As for the bottom row, we have
$$
\bigsqcup_{p=1}^{n-1} (-1)^{p+q} \binom{[n+p-j]}{n-1} \binom{[n]}{q-p} =
\bigsqcup_{k=0}^{n} (-1)^{k} \binom{[n+q-k-j]}{n-1} \binom{[n]}{k} = \begin{cases} \si 0 0 & q=j-1 \\ \u \emptyset &
\text{otherwise} \end{cases},
$$
and we obtain a sijection to $(-1)^{n+j} \det( \u {\p T}^j)$.

%We set $\u {\p S} = [(-1)^{i+j} \binom{[n]}{j-i}]_{2 \le i, j \le n}$ and $\u {\p S}' = [\binom{[n+j-i-1]}{j-i}]_{2 \le i, j \le n}$, and use Problem~\ref{prob:binom}
%to see that there is a sijection
%$$
%\bigsqcup_{k=2}^n (\u {\p S})_{i,k} \times (\u {\p S}')_{k,j} \Longrightarrow \begin{cases} [0,0] & i=j \\ \emptyset & \text{otherwise} \end{cases}.
%$$
%Using Problem~\ref{prob:detproduct}, it follows that there is a sijection
%$$
%\DPP_{n-1} \times \ASM_{n,j} \Longrightarrow \det \left( \u {\p S} \cdot \u {\p P} \cdot \u {\p S}' \sqcup \u {\p S} \cdot \u {\p Q}^j \cdot \u {\p S}' \right) \sqcup
%- \det \left( \u {\p P} \right).
%$$
%Employing Problem~\ref{prob:binom} several times, a sijection
%$$
%\det \left( \u {\p S} \cdot \u {\p P} \cdot \u {\p S}' \sqcup \u {\p S} \cdot \u {\p Q}^j \cdot \u {\p S}' \right) \sqcup
%- \det \left( \u {\p P} \right) \Longrightarrow \det \left( \u { \p R }^j \right)
%$$
%where $\u { \p R }^j = [\u R_{p,q}]_{1 \le p, q \le n-1}$ with
%$$
%\u R_{p,q} = \begin{cases} \binom{[i+j]}{j-1} & p < n-1 \\ \binom{[n+q-j]}{n-1} & p=n-1, j>1 \\ 0 & p=n-1, j=1 \end{cases} \sqcup \begin{cases} [0,0] & i<n-1 \text{ or } k=1 \\ \emptyset & \text{otherwise} \end{cases}.
%$$
%From Lindstr\"om-Gessel-Viennot it follows that there is a sijection
%$$
%\det \left( \u { \p R }^j \right)  \Longrightarrow \DPP_{n,j}.
%$$
\end{proof}

\section{Concluding remarks}

In this paper, we present the first bijective proof of the enumeration formula for alternating sign matrices. The bijection is by no means simple; the construction is based on  \cite{PartI}, and combined the two papers have about 40 pages, with the technical constructions taking about 20 pages. We also needed more than 2000 lines to produce a working Python code. However, note that the first proof of the ASM theorem by Zeilberger was 84 pages long. We certainly hope that our proof will be simplified and shortened in the future.

On the other hand, this successful translation of a computational proof into a bijective proof also raises questions as to the relation between these two types of proofs. Under what circumstances is such a translation possible? If there exists a simple explicit bijection between, say, ASMs and DPPs, can it be converted into a (simple) computational proof? 

\section*{Acknowledgments}

We are grateful to Matija Pretnar, Alex Simpson and Doron Zeilberger for interesting conversations and helpful suggestions.

%\bibliographystyle{alpha}
%\bibliography{asmpp}

\end{document}